\documentclass[11pt]{article}

\PassOptionsToPackage{usenames,dvipsnames}{xcolor}
\usepackage{xcolor}

\usepackage[margin=1.25in]{geometry} 
\usepackage{graphicx} 
\usepackage{tikz}
\usetikzlibrary{decorations.markings}
\usepackage{booktabs} 
\usepackage{array} 
\usepackage{paralist} 
\usepackage{verbatim} 
\usepackage{mathrsfs} 
\usepackage{amssymb}
\usepackage{amsthm}
\usepackage{amsmath,amsfonts,amssymb}
\usepackage{esint}
\usepackage{graphics}
\usepackage{enumerate}
\usepackage{mathtools}
\usepackage{xfrac}
\usepackage{bbm}
\usepackage{subcaption}
\usepackage{stmaryrd} 

\let\oldsquare\square 

\usepackage{mathabx}

\renewcommand{\square}{\oldsquare}

\usepackage{tocloft}

\newcommand{\ga}{\gamma}

\newcommand{\la}{\lambda}
\newcommand{\NN}{\mathbb{N}}
\newcommand{\CC}{\mathbb{C}}

\newcommand{\sh}{\mathcal{H}}

\numberwithin{equation}{section}

\newtheorem{theorem}{Theorem}[section]

\newtheorem{corollary}[theorem]{Corollary}
\newtheorem{proposition}[theorem]{Proposition}
\newtheorem{lemma}[theorem]{Lemma}
\theoremstyle{definition}
\newtheorem{definition}[theorem]{Definition}

\newtheorem{remark}[theorem]{Remark}

\let\originalleft\left
\let\originalright\right
\renewcommand{\left}{\mathopen{}\mathclose\bgroup\originalleft}
\renewcommand{\right}{\aftergroup\egroup\originalright}

\newcommand{\vertiii}{\vert\kern-0.3ex\vert\kern-0.25ex\vert}

\newcommand*{\Z}{\ensuremath{\mathbb{Z}}}
\newcommand*{\R}{\ensuremath{\mathbb{R}}}

\newcommand*{\Rd}{\ensuremath{\mathbb{R}^d}}

\renewcommand*{\tilde}{\widetilde}

\DeclareSymbolFont{boldoperators}{OT1}{cmr}{bx}{n}
\SetSymbolFont{boldoperators}{bold}{OT1}{cmr}{bx}{n}

\renewcommand{\P}{\mathbb{P}}

\newcommand{\bE}{\mathbf{E}}
\newcommand{\bF}{\mathbf{F}}
\newcommand{\indc}{1}

\newcommand{\bH}{\mathbf{H}}

\DeclareMathOperator*{\argmin}{argmin}

\newcommand{\RR}{\mathbb{R}}
\newcommand{\de}{\delta}
\newcommand{\ENZ}{\Omega \setminus \overline{D}}
\newcommand{\e}{\varepsilon}

\makeatletter
\newcommand\avsuminner[2]{%
	{\sbox0{$\m@th#1\sum$}%
		\vphantom{\usebox0}%
		\ooalign{%
			\hidewidth
			\smash{\,\rule[.23em]{8.8pt}{1.1pt} \relax}%
			\hidewidth\cr
			$\m@th#1\sum$\cr
		}%
	}%
}
\makeatother

\makeatletter
\newcommand\avsuminnerr[2]{%
	{\sbox0{$\m@th#1\sum$}%
		\vphantom{\usebox0}%
		\ooalign{%
			\hidewidth
			\smash{\,\rule[.23em]{6pt}{0.7pt} \relax}%
			\hidewidth\cr
			$\m@th#1\sum$\cr
		}%
	}%
}
\makeatother

\def\Xint#1{\mathchoice
	{\XXint\displaystyle\textstyle{#1}}%
	{\XXint\textstyle\scriptstyle{#1}}%
	{\XXint\scriptstyle\scriptscriptstyle{#1}}%
	{\XXint\scriptscriptstyle\scriptscriptstyle{#1}}%
	\!\int}
\def\XXint#1#2#3{{\setbox0=\hbox{$#1{#2#3}{\int}$}
		\vcenter{\hbox{$#2#3$}}\kern-.5\wd0}}

\def\fint{\Xint-}

\newcommand{\subjclass}[2][]{%
	\par\smallskip
	\noindent\textbf{Mathematics Subject Classification (2020):}
	\if\relax\detokenize{#1}\relax
	#2%
	\else
	#1, #2%
	\fi
	\par\smallskip
}

\newcommand{\keywords}[1]{%
	\par\smallskip
	\noindent\textbf{Keywords:} #1
	\par\smallskip
}

\usepackage[colorlinks=true, pdfstartview=FitV,
linkcolor=blue, citecolor=blue, urlcolor=blue]{hyperref}

\usepackage{titlesec}

\newcommand{\addperiod}[1]{#1.}
\titleformat{\section}
{\normalfont\Large}{\thesection.}{0.5em}{}
\titleformat*{\subsection}{\normalfont\large}
\titleformat{\subsubsection}[runin]
{\bfseries}
{\thesubsubsection.}
{0.5em}
{\addperiod}
\titleformat*{\subsubsection}{\bfseries}
\titleformat*{\paragraph}{\bfseries}
\titleformat*{\subparagraph}{\large\bfseries}

\title{\bf \Large Analytic spectral perturbation theory for a high-contrast Maxwell operator}

\author{Robert V. Kohn
	\thanks{(Deceased) Courant Institute of Mathematical Sciences, New York University.
		{\footnotesize \href{mailto:kohn@cims.nyu.edu}{kohn@cims.nyu.edu}. }
	}
	\and 
	Raghavendra Venkatraman
	\thanks{Department of Mathematics, The University of Utah.
		{\footnotesize \href{mailto:raghav@math.utah.edu}{raghav@math.utah.edu}. \\
			The authors warmly thank Nader Engheta for bringing this problem to our attention, and for many
			helpful discussions. Both authors gratefully acknowledge support from the Simons Foundation through its
			Collaboration on Extreme Wave Phenomena (grant 733694). RVK also gratefully acknowledges support
			from the National Science Foundation (grant DMS-2009746). RV also gratefully acknowledges support from the National Science Foundation (grant DMS-2407592)}
	}
}

\date{\today}

\usepackage[nottoc,notlot,notlof]{tocbibind}

\begin{document}
	
	\maketitle
	
\begin{abstract}
	We study analytic spectral perturbation theory for the time-harmonic Maxwell operator in a perfectly electrically conducting cavity containing a high-contrast core--shell structure. The dielectric permittivity equals $1$ in a bounded inclusion and a small complex parameter $\delta$ in the surrounding shell. The limit $\delta \to 0$ corresponds to an infinite-contrast regime and leads to a degenerate Maxwell system. Despite this degeneracy, we develop a detailed spectral theory for the limiting problem for general Lipschitz inclusions and shells.
	
	Using a novel operator-theoretic reformulation, we prove complex-analytic dependence of the spectrum on $\delta$ in a neighborhood of $\delta = 0$. When the inclusion is a ball, we analyze the asymptotic expansion of eigenvalues and identify conditions under which the leading-order term is independent of the geometry of the surrounding shell. We also construct examples of resonances for which the leading-order asymptotics depend sensitively on the shell geometry, even in this symmetric setting. These results clarify the mechanisms underlying geometry-invariance of resonances in high-contrast Maxwell systems and explain their robustness under small complex perturbations.
\end{abstract}


	\setcounter{tocdepth}{2}  
	\tableofcontents
\section{Introduction} \label{sec:intro}
\subsection{Motivation and Background}
In this paper, we analyze a family of eigenvalue problems for the Maxwell system in a bounded, smooth domain~$\Omega \subset \R^3,$ where the dielectric permittivity~$\e = \e_\delta(x)$ nearly vanishes in part of the domain~$\Omega.$ To explain the set-up, we let~$D \subset \Omega$ be a smooth sub-domain, and study the cavity eigenproblem for the Maxwell system with a perfectly electrically conducting (PEC) boundary: 
\begin{equation} \label{e.introPDEfront}
	\begin{dcases}
		\nabla \times \nabla \times \bE (x)&= \lambda \e_\delta (x) \bE(x)\, \quad x\mbox{ in } \Omega\\
		\bE(x) \times \nu_\Omega(x) &= 0, \qquad \qquad x \mbox{ on } \partial \Omega\,;
	\end{dcases}
\end{equation}

We are interested \textit{high-contrast} complex dielectric permittivities of the form~$\e_\delta(x) = 1$ if~$x \in D,$ the dielectric inclusion, and~$\e_\delta(x) = \delta$ if~$x \in \ENZ,$ where~$\de \in \mathbb{C}$ is a complex number with~$|\delta|\ll 1$. In the above,~$\nu_\Omega$ denotes the outward unit normal to~$\Omega,$ so that the boundary conditions in the eigenvalue problem above require that~$\bE$ is normal to the outer boundary~$\partial\Omega.$ A first glance at the eigenvalue problem~\eqref{e.introPDEfront} reveals that it appears somewhat singular: since~$\e_\de$ is not bounded away from zero, the equation gives only \textit{fading} control on the divergence of~$\bE$ on~$\Omega$, as it only asserts that~$\mathrm{div}\, (\e_\de \bE),$ being the divergence of a curl, vanishes.  

We are interested in the dependence of the spectrum on~$\de$ in a neighborhood of~$\de=0.$ As we explain below, one motivation to study~\eqref{e.introPDEfront} in the limit~$\de \to 0$ is that for \textit{certain} choices of~$D$, the eigenvalue problem~\eqref{e.introPDEfront} has eigenvalues, which, to leading order in~$\de$, \textit{are independent of the geometry of~$\Omega.$} Understanding the spectral perturbation theory of the family of problems~\eqref{e.introPDEfront}, and explaining the \textit{character} and \textit{robustness} of the asymptotic geometry invariance of eigenvalues of this Maxwell system are the primary goals of this paper. The principal contributions of our paper are, roughly speaking, \textit{complex analytic dependence} of the spectrum  on~$\de$ in a neighborhood of~$\de = 0,$ and providing explicit PDE characterizations for each term of the Taylor expansion. As is always the case in spectral perturbation theory, care is required while discussing perturbation theory about a multiple eigenvalue of the limiting problem, and our analysis does cover this case.

\begin{figure}[htbp]
	\centering 
	\includegraphics[scale=.9]{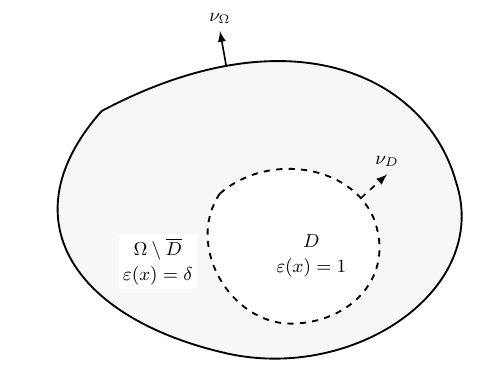}
		\caption{A core-shell resonator with~$\e(x) \equiv \delta$ in the shell~$\Omega \setminus \overline{D}$ and~$\e(x) \equiv 1$ in the inclusion D. Also marked are outward unit normals to both the inclusion and the shell. }
		\label{f.pic}
\end{figure}

Let us begin with some physical motivation and background for our interest in this problem, referring the mathematically inclined reader to Section~\ref{ss.mainresults} where we precisely lay out our assumptions and describe our main results.  The problem considered here was introduced in a 2016 Nature Communications paper by Liberal et.al.~\cite{Engheta}. That paper belongs to a broader, emerging literature (see, for instance,~\cite{silveirinha2006tunneling,silveirinha2007design,silveirinha2007theory,liberal2017rise,niu2018epsilon} and references therein) in applied physics (photonics) which explores the striking scattering and resonance properties of devices that operate at a frequency at which one of the components of the device has zero dielectric permittivity. One such property, the primary motivating application in~\cite{Engheta} as well as in the present paper, is that under certain settings, a perfectly electrically conducting (PEC) core-shell cavity consisting of a lossless dielectric core surrounded by an ``epsilon-near-zero'' (ENZ) shell has a \textit{resonance} at the ENZ frequency that is  \textit{independent of the geometry of the shell in the lossless limit}, see Fig~\ref{f.pic}. The present paper, along with our previous paper~\cite{KV2}, succeed in putting the discussions of~\cite{Engheta} on a sound mathematical foundation. In yet another recent paper~\cite{KV1}, we study the unconventional scattering properties of an ENZ scatterer with a dielectric inclusion. 

There are two settings under which the authors of~\cite{Engheta} propose geometry-invariant resonant cavities in the lossless limit: a transverse magnetic (TM) setting that leads to the study of a two-dimensional problem that we analyzed in~\cite{KV2}, and a fully three-dimensional eigenproblem for Maxwell's system that we analyze here. Before describing the latter, let us begin by describing the geometry-invariant character of the resonances in the (simpler) TM setting to aid in intuition. 

 Mathematically this corresponds to studying time-harmonic Maxwell's equations in a cylindrical domain along the~$e_3-$axis with cross-sectional shape~$\Omega,$ which hosts a cylindrical dielectric inclusion with cross section~$D,$  in which the magnetic field has the form~$H(x_1,x_2,x_3) = (0,0,u(x_1,x_2))$ for an unknown complex-scalar field~$u.$ The eigenproblem for Maxwell's equations with a perfectly electrically conducting boundary condition on~$\partial \Omega$ reduce to a Neumann eigenproblem
 \begin{equation} \label{e.2Dproblem}
 	\begin{dcases}
 	&	- \mathrm{div}\, \e_\de^{-1}(x,\omega) \nabla u = \omega^2 u, \quad \mbox{ in } \Omega\\
 	&	\partial_{\nu_\Omega} u = 0, \qquad \qquad \quad \qquad \mbox{ on } \partial \Omega\,,
 	\end{dcases}
 \end{equation}
where~$\omega$ is the operating frequency of the device, and~$\e_\de$ denotes dielectric permittivity (as with~\cite{Engheta}, we assume that we are working with non-magnetic materials, and so, we have set the magnetic permeability~$\mu \equiv 1$ in both the shell and the dielectric inclusion). It is important to note that~\eqref{e.2Dproblem} is a nonlinear eigenvalue problem since the resonance frequency~$\omega$ also appears in the coefficient~$\e_\de(x,\omega):$ this is the phenomenon of dispersion. Ignoring dispersion for the moment and studying the resulting \textit{linear} Neumann problem~\eqref{e.2Dproblem} where we set~$\e_\de(x) \equiv 1$ for~$x \in D,$ and~$\e_\de(x) \approx 0$ for~$x \in \ENZ,$ it is already evident what is so special about the~$\de \to 0$ limit: the vector field~$\e^{-1} \nabla u$ can avoid being large provided~$\nabla u \approx 0$ in the region~$\ENZ$ where~$\e(x) = \de \in \CC, \de \approx 0.$ In the limit~$\de \to 0^+,$ the magnetic field~$u$ is \textit{constant}, and so the problem of solving the PDE in~\eqref{e.2Dproblem} reduces to determining this constant. It is this insight that is at the heart of the geometry invariance of such a TM resonance: for a given inclusion shape~$D,$ the value of this constant magnetic field in the~$\de \to 0$ limit turns out to depend only on the \textit{area} of the ENZ region~$\ENZ$, not its shape/geometry. 

The ENZ limit is an idealization, for two reasons, both related to the phenomenon of \textit{dispersion} (see~\cite[Section 82]{LL}). In some detail, the (complex) dielectric permittivity is a function of frequencies~$\omega \in \CC$ with certain structural properties reflecting physical considerations such as causality.  Real materials always have losses that result in a small non-zero imaginary part of the dielectric permittivity~$\e_\de(x)$. Furthermore, the real part of ~$\e$ is exactly zero only at isolated frequencies (referred to in~\cite{Engheta} as ENZ frequencies), and so if one desires a broadband device that operates at frequencies \textit{near}, rather than exactly at the ENZ frequency, then even the real part of~$\e$ is merely small, and not zero. 

Due to the presence of the reciprocal permittivity~$\frac1{\e_\de} = \frac1{\delta}$ in the ENZ region, the PDE problem~\eqref{e.2Dproblem} appears to be rather singular on first sight. Nonetheless, in our previous paper~\cite{KV2}, we show, roughly speaking, that the spectrum of~\eqref{e.2Dproblem} depends in a complex analytic manner on~$\de$ in a complex neighborhood of~$\de=0.$ More precisely, we analyze~\eqref{e.2Dproblem} in two steps: first, ignoring dispersion, we use an argument with the flavor of spectral perturbation theory, and demonstrate the analytic dependence of the spectrum near~$\de = 0$, for the \textit{linear} eigenvalue problem 
\begin{equation*}
	\begin{dcases}
		&- \mathrm{div} \bigl( \indc_D + \delta^{-1} \indc_{\ENZ}\bigr) \nabla u = \lambda(\delta) u, \quad \mbox{ in } \Omega\\
	&	\partial_{\nu_\Omega} u = 0 \qquad \qquad \qquad \qquad \qquad \qquad \mbox{ on } \partial \Omega\,. 
	\end{dcases}
\end{equation*}
Then, perturbing about the limit of zero loss, for a fairly general dispersion relation~$\e = \e(\omega)$ that includes the Lorentz and Drude models, we use the implicit function theorem to argue that the spectrum of the nonlinear eigenvalue problem~\eqref{e.2Dproblem} depends analytically on the (small) loss parameter. 

Before ending our discussion on the TM problem, let us emphasize that even though the \textit{leading order} value of the eigenvalue, corresponding to the~$\de = 0$ limit, depends only on the area of the ENZ region~$\ENZ,$ the corrections to~$\lambda(\delta)$ when~$\de \neq 0$, \textit{do depend on the geometry of~$\Omega$,} and in fact, have explicit PDE characterizations. 
In fact, in~\cite{KV2} we also analyze a planar shape optimization problem that maximizes the leading order correction to~$\lambda(\de),$ given by $ \lambda'(\delta)|_{\delta = 0}.$ Elementary perturbation theory reveals that this quantity is always non-positive and depends on the shape of~$\Omega$. When~$\de$ has a non-zero imaginary part corresponding to losses, this quantity dictates, to leading order, how fast this ``resonance'' decays in time. In light of this, in~\cite{KV2} we also address the question: for a given shape of the inclusion~$D$, what is the optimal shape of the ENZ shell~$\Omega$, that results in a resonator whose resonance decays the slowest in time. 

We turn now to the setting of our present paper. This corresponds to the second type of geometry-invariant resonant cavities introduced by the authors in~\cite{Engheta}. To describe these, given~$\Omega \subset \R^3$ (the shell) and the dielectric inclusion~$D,$ contained in the interior of~$\Omega$ (see Fig.~\ref{f.pic}), we are interested in an eigenvalue~$\omega$ and nonzero electric and magnetic fields~$E$ and~$H,$ respectively, which satisfy 
\begin{equation} \label{e.maxwell-intro}
	\begin{dcases}
		\nabla \times E = - i \omega H, \quad \mbox{ in } \Omega\\
		\nabla \times H = i \omega \e E, \quad \mbox{ in } \Omega \\
		\nu_\Omega \times E = 0, \quad \quad \mbox{ on } \partial \Omega\,,
	\end{dcases}
\end{equation} 
where, once again, the dielectric permittivity~$\e(x) = 1$ for~$x \in D$, and~$\e(x) = \delta$ for~$x \in \ENZ,$ where~$\delta \in \CC$ with~$|\delta|\ll 1.$ Thus, the dielectric permittivity depends on the small, complex-valued parameter~$\de,$ and so we often write~$\e_\de$ to indicate this dependence. Ignoring dispersion, once again, when~$\e(x) \approx 0$ then the second equation in~\eqref{e.maxwell-intro} asserts that the magnetic field~$H$ is very nearly curl-free in the ENZ region. Observe that, in contrast with the TM setting discussed above, in the limit~$\delta \to 0,$ this curl-free (i.e., gradient) constraint on~$H$ results in an infinite-dimensional set of possibilities for the magnetic field in the ENZ region, and furthermore, it is not clear if such resonances, should they exist, can be independent of the geometry of~$\Omega.$ Furthermore, once again, from the second Maxwell equation, since the vector field~$\e E$ is (a constant multiple of) a curl, it must be divergence-free: this entails that in the~$\delta \to 0$ limit, the electric field~$E$, which is required to satisfy the \textit{normal} PEC boundary condition~$\nu_\Omega \times E = 0$ at the outer boundary~$\partial \Omega,$ is also required to satisfy the \textit{tangential} boundary condition~$E \cdot \nu_D= 0$ at the transmission interface~$\partial D.$ This change of the electric field, too, from tangential at~$\partial D$ to normal at~$\partial \Omega,$ for the solution to a PDE, also appears to be rather sensitive to the geometry of the domain~$\Omega.$ 

In spite of the geometric constraints explained in the preceding paragraph, remarkably, the authors in~\cite{Engheta} discover solutions to~\eqref{e.maxwell-intro} \textit{in the~$\e \to 0$ limit,} for specific choices of the inner dielectric inclusion~$D.$  For instance, they reason, through separation of variables, that when the inclusion~$D$ is a ball,  there exist solutions to this system which have the property that when~$\e = 0$, the magnetic field~$H$ is not merely curl-free in~$\ENZ$; it \textit{vanishes identically} in the ENZ region~$\ENZ,$ while being non-vanishing in the inclusion~$D$. For such a magnetic field which identically vanishes in~$\ENZ,$ the electric field is naturally curl-free (to leading order), and hence, the electric field~$E,$ which is individually divergence free in each of the inclusion and the ENZ shell, is the gradient of a harmonic function. For this reason, they refer to such resonances as \textit{electrostatic resonances}. 

Summarizing, the electric and magnetic fields corresponding to such an electrostatic resonance must satisfy an over-determined problem in the inclusion~$D$: 
\begin{equation} \label{e.overdet}
	\begin{dcases}
			\nabla \times E = - i \omega H, \quad \mbox{ in } D\\
		\nabla \times H = i \omega E, \ \quad \mbox{ in } D \\
		\nu_\Omega \cdot E = 0, \ \ \ \ \quad \quad \mbox{ on } \partial D\\
		H = 0, \qquad \qquad \quad \mbox{ on } \partial D\,.
	\end{dcases}
\end{equation}
Given such a solution to the overdetermined eigenvalue problem in~\eqref{e.overdet} (which, evidently is independent of the geometry of~$\Omega$) it is not hard to extend this to a solution to~\eqref{e.maxwell-intro} on the whole set~$\Omega,$ with the same frequency (see Section~\ref{sec:electro}).  

The foregoing discussion on the Maxwell eigenproblem ignores dispersion. However, as explained in~\cite{KV2}, it is straight-forward to include the effects of dispersion through the implicit function theorem, once we have shown the analyticity of the spectrum of the \textit{linear} eigenvalue problem in~\eqref{e.maxwell-intro} with~$\e = \e_\de$ independent of~$\omega.$ For this reason, we exclusively focus on this latter step in the rest of the paper, and do not consider dispersion. 

\subsection{Problem set-up and main results}
\label{ss.mainresults}
We let $\Omega \subset \R^3$ denote a connected, simply-connected Lipschitz domain, and let $D$ be a connected Lipschitz subdomain which is such that $\overline{D} \subset \Omega.$ For any~$\de \in \CC$, we let $\e_\de$ be defined via 
\begin{equation} \label{e.ENZmedium}
	\e_\de(x) := \left\{
	\begin{array}{cc}
		1 & \mbox{ if } x \in D\\
		\de &  \mbox{ if } x \in \ENZ\,,
	\end{array}
	\right.
\end{equation}
or equivalently,
\begin{equation} \label{e.epsdeltaexp}
	\e_\de = \indc_D + \de \indc_{\ENZ}\,. 
\end{equation}
For any $\de \in \mathbb{C} \setminus \{0\}$ which is such that $|\de| \ll 1,$  we are interested in eigenvalues $\la_\de \in \mathbb{C} \setminus \{0\}$ and associated eigenfunctions $\bE \in \mathrm{H}_0(\mathrm{curl}\,, \Omega)$ of the Maxwell's equation 
\begin{equation} \label{e.maxwell-delta}
	\int_\Omega \nabla \times \bE \cdot \overline{\nabla \times \bF}\,dx - \lambda_\de \int_\Omega \e_\de \bE \cdot \overline{\bF}\,dx = 0\,, \mbox{ for all } \bF \in \mathrm{H}_0(\mathrm{curl}\,, \Omega)\,.
\end{equation}
Here, we recall that 
\begin{equation*}
	H_0(\mathrm{curl},\Omega) := \{f \in L^2(\Omega,\R^3) : \nabla \times f \in L^2(\Omega,\R^3) \mbox{ and } \nu_\Omega \times f = 0 \mbox{ on } \partial \Omega\}\,. 
\end{equation*}
This generalized eigenvalue problem has strong formulation given by
\begin{equation}
	\label{e.strongform}
	\begin{dcases}
		\nabla \times \nabla \times \bE_{\de} &= \la_{\de} \e_\de \bE_\de, \quad \quad \mbox{ in } \Omega\\
		\bE_\de \times \nu_\Omega &= 0  \quad \quad \quad \quad \mbox{ on } \partial \Omega\,. 
	\end{dcases}
\end{equation}
The ``limiting problem'' associated with~$\de=0$ is the generalized eigenvalue problem 
\begin{equation}
	\label{e.strongform-lim}
	\begin{dcases}
		\nabla \times \nabla \times \bE_{0} &= \la_{0} \indc_D \bE_0, \quad \quad \mbox{ in } \Omega\\
		\bE_0 \times \nu_\Omega &= 0  \quad \quad \quad \quad \mbox{ on } \partial \Omega\,. 
	\end{dcases}
\end{equation}
The problem~\eqref{e.strongform} is a family of generalized eigenvalue problems of the form~$Av = \lambda B_\de v,$ where, informally,~$A$ is the Maxwell operator with the homogeneous tangential component boundary conditions, and~$B$ is the operator that multiplies a vector field by~$\e_\delta \in L^\infty(\Omega)$ (which nearly vanishes on the shell). Due to the degenerate nature of the permittivity~$\e_\de$, the equations~\eqref{e.strongform}-\eqref{e.strongform-lim} do not provide control on~$\mathrm{div}\, \bE_\de$ in the ENZ region~$\ENZ$. In particular, the standard functional analytic framework that applies the Lax-Milgram theorem to a suitable \textit{coercive} bilinear form to study the spectral theory of the Maxwell operator becomes unavailable; see~\cite{Cos,Monk}. 

In light of the vectorial nature of the problem, it also is inconvenient to study the problem ``region-by-region'' as we did in our analysis of the TM problem in~\cite{KV2}. 
In order to circumvent these difficulties, we study the eigenvalue problems~\eqref{e.strongform}-\eqref{e.strongform-lim} through a suitable factorization of the equation; this results in a family of compact, self-adjoint operators on a suitable Hilbert-space that depend analytically on a small parameter~$\de.$ We show that this family of compact self-adjoint operators is such that Rellich's perturbation theorem (recalled in the Appendix) applies. 

Let us explain our strategy in more detail in the context of a simpler, scalar problem which shares several of the features of the vectorial problem; we provide a detailed analysis of this scalar problem in Section~\ref{s.scalar} below. In that section, we study perturbation theory for the family of Neumann eigenvalue problems 
\begin{equation} \label{e.PDEscalar-intro}
	\begin{dcases}
		-\Delta u_\delta &= \lambda_\delta \e_\delta u_\delta\,  \quad \mbox{ in } \Omega\\
		\partial_{\nu_\Omega} u_\de &= 0, \qquad \quad \mbox{ on }\partial \Omega\,,
	\end{dcases}
\end{equation}
about their~$\de=0$ limiting problem:
\begin{equation} \label{e.PDEscalar-0-intro}
	\begin{dcases}
		-\Delta u &= \lambda \indc_D u\,  \quad \mbox{ in } \Omega\\
		\partial_{\nu_\Omega} u &= 0, \qquad \quad \mbox{ on }\partial \Omega\,. 
	\end{dcases}
\end{equation}
In order to study describe our approach to the spectral theory of~\eqref{e.PDEscalar-0-intro} we recast it as follows. Writing~$\mathring{L}^2(\Omega)$ for mean zero functions in~$\Omega$, consider the operator~$K_0$ defined via 
\begin{equation} \label{e.K0defintro}
	K_0 := (-\Delta)^{-\sfrac12} \indc_D (-\Delta)^{-\sfrac12}\,,
\end{equation}
where~$\indc_D$ is the multiplication operator by the bounded function~$\indc_D$, the indicator function of the inclusion domain~$D$, defined via~$\indc_D(x) = 1$ when~$x \in D, \indc_D(x) = 0$ when~$x\not\in D.$ Here,~$(-\Delta)^{-1}: \mathring{L}^2(\Omega) \to \mathring{L}^2(\Omega)$ is the unique mean-zero solution  to the Neumann problem 
\begin{equation*}
	\begin{dcases}
		-\Delta u & = f \quad \mbox{ in } \Omega\\
		\partial_{\nu_\Omega} u &= 0, \quad \mbox{ on } \partial \Omega\,,
	\end{dcases}
\end{equation*} 
given data~$f\in \mathring{L}^2(\Omega);$ and we set~$u = (-\Delta)^{-1} f.$ With this notation, we define fractional powers of this positive operator in the usual, spectral manner; see Section~\ref{s.scalar} for the precise definition.  One must then extend the operator~$K_0$ to the full space~$L^2(\Omega)$ through projections; this is done in~\eqref{e.projectiondefK_0}. Eigenvalues for~\eqref{e.PDEscalar-0-intro} are obtained by taking inverses of the nonzero eigenvalues of the compact operator~$K_0.$ The associated eigenfunctions of~$K_0$, upon adding a suitable constant (which the Neumann boundary condition at~$\partial \Omega$ does not see), yield  eigenfunctions to the limiting PDE problem~\eqref{e.PDEscalar-0-intro}. In this way, the generalized eigenvalue problem~\eqref{e.PDEscalar-0-intro} can be put into a framework to which Fredholm theory applies, see Theorem~\ref{t.limscalar} for a precise formulation. 

Still in the scalar case, the viewpoint on the~$\de \neq 0$ family of problems, in~\eqref{e.PDEscalar-intro} is similar and proceeds indirectly through a compact hermitian operator~$K_\delta$, introduced in~\eqref{e.Kdelta}. Care must be taken, however, in defining the projections that go into defining~$K_\delta$, since these, too, depend on~$\de,$ but \textit{analytically}. The upshot of this discussion is that~$\{K_\de\}_{\de}$ is a family of compact self-adjoint operators on~$L^2(\Omega)$ that depend in a complex analytic manner on~$\de$ in a neighborhood of~$\de=0,$ with the property that when restricted to~$\de \in \RR$ near~$\de=0,$ the operators~$K_\de$ are also self-adjoint. Therefore, we can appeal to a deep theorem of F. Rellich, recalled in Appendix~\ref{sec:rellich}. Starting from an eigenfunction of~$K_\de$ (which, by the foregoing analytic perturbation theory argument, depends analytically on~$\de$), addition of a constant, which too, depends analytically on~$\de$, yields an eigenfunction for the generalized eigenvalue problem~\eqref{e.PDEscalar-intro}. Our result for the complex analytic dependence of the spectrum of the scalar problem~\eqref{e.PDEscalar-intro} appears in Theorem~\ref{t.scalar}. 

Our analysis of the Maxwell eigenvalue problems~\eqref{e.strongform}-\eqref{e.strongform-lim} appears in Section~\ref{sec:vec}. While the overall structure of the argument follows the scalar case, there are a number of complications that arise due to the fact that the kernel of the~$\nabla \times\nabla \times $ operator with homogeneous tangential (PEC) boundary conditions is infinite dimensional, and consists of arbitrary gradients of~$H^1$ functions that are constant on~$\partial \Omega.$ In a little detail, by analogy with~\eqref{e.K0defintro}, it is natural to study the analytic family of operators given by 
\begin{equation*}
	\mathbb{K}_\delta := \bigl( \nabla \times \nabla \times )^{-\sfrac12} \e_\delta \bigl( \nabla \times \nabla \times )^{-\sfrac12}\,,
\end{equation*}
defined on~$L^2$ vector fields, where, once again,~$(\nabla \times \nabla \times)^{-\sfrac12}$ is defined via the spectral theorem applied to the positive solution operator that inverts~$\nabla \times \nabla \times u =f$ with PEC boundary conditions. Evidently, this solution operator acts on \textit{divergence-free} fields~$f$. However, the presence of the discontinuous \textit{degenerate elliptic} coefficient field~$\e_\de$ necessitates a family of~$\delta-$dependent Helmholtz-Hodge projections which project a given vector field~$F$ onto a space of vector fields such that~$\e_\de F$ is divergence-free. In turn, this necessitates studying the solutions to an elliptic equation whose coefficients are nearly zero (and complex!) on part of the domain. At the heart of the matter is Lemma~\ref{l.analyticity-poisson-highcontrast} which proves analyticity of a family of Helmholtz-Hodge decomposition with degenerate metrics; its proof relies on an infinite order expansion procedure. We highlight that since our arguments avoid the use of layer potential theoretic arguments, it works for Lipschitz domains with Lipschitz inclusions. 

With this technical lemma in place, the analysis of the limiting eigenproblem and its spectral theory appears in Theorem~\ref{t.maxwelllim-spectrum}, and its spectral perturbation theory appears in Section~\ref{s.perturb-vector}, culminating in Theorem~\ref{t.vectormain}. These theorems demonstrate that, roughly speaking, the eigenvalues of the Maxwell eigenproblem~\eqref{e.introPDEfront} depend \textit{complex analytically} on~$\de$ in a neighborhood of~$\de = 0$-- an analytic explanation of the robustness of the effects associated with the permittivity being exactly zero in~$\ENZ.$ In~\cite{Engheta}, the authors investigate such robustness mainly through numerical simulation. 

Having studied the spectral theory of the Maxwell generalized eigenproblem~\eqref{e.introPDEfront}, in Section~\ref{sec:electro} we finally address the existence of \textit{electrostatic} resonances to the eigenvalue problem~\eqref{e.introPDEfront}-- namely, eigenfunctions~$\bE$ which satisfy~$\nabla \times \bE = O(\de)$ in the~$\ENZ$ region, so that the associated magnetic fields \textit{vanish} to leading order in the ENZ region~$\ENZ.$ We show in this section that 
\begin{enumerate}
	\item When the inclusion~$D$ is a ball, there \textit{do exist} infinitely many electrostatic eigenvalues. As explained earlier, these are solutions to the overdetermined problem for Maxwell's equations in~\eqref{e.overdet}. Thus, in this case, in agreement with~\cite{Engheta} we produce an infinite family of high multiplicity eigenfunctions of~\eqref{e.introPDEfront} whose~$\de \to 0$ limit is \textit{independent of the shape of}~$\Omega.$ This is the content of Theorem~\ref{t.electro}.
	\item Still when~$D$ is a ball, and given~$\Omega$ containing~$D$, we show the existence of \textit{infinitely} many \textit{non-electrostatic} resonances, i.e., eigenfunctions whose associated magnetic field~$H$ is \textit{non-vanishing} to leading order in~$\ENZ,$ but are merely \textit{curl-free}. This is the content of Theorem~\ref{p.notelectrostatic}.
\end{enumerate}

As explained above, our approach relies on an abstract operator theoretic formulation of~\eqref{e.introPDEfront} which allows for a general, quantitative theory that does not depend on special symmetries of the geometry in which separation of variables is available. This offers an understanding of~\cite{Engheta} as to (a) the character of the geometry invariance of resonances, and (b) the robustness of these resonances to losses; the authors in~\cite{Engheta} offered a different type of understanding of (b)-- one that was based on numerical simulation. 

Let us point out that the overdetermined problem~\eqref{e.overdet} bears resemblance with the celebrated overdetermined problem of Serrin~\cite{serrin} (see also \cite{weinberger}): if~$\Omega \subset \Rd$ is a bounded, open, smooth, connected domain in~$\Rd,$ and~$u$ solves the overdetermined problem 
\begin{equation}
	\begin{dcases}
		\Delta u &= -1 \quad \mbox{ in } \Omega\\
		u &= 0 \quad \quad \mbox{ on } \partial \Omega\\
		\frac{\partial u}{\partial \nu}  &= \mbox{ constant } \quad \quad \mbox{ on } \partial \Omega\,,
	\end{dcases}
\end{equation}
then, in fact,~$\Omega$ must be a ball, and~$u$ must be a radial function about the center of the ball. Formulating Maxwell's equations in terms of differential forms, note that the overdetermined probem~\eqref{e.overdet} in the inclusion~$D$ (with~$\e = 1$) writes as: we seek~$\lambda > 0$ and a square integrable one form~$\bE$ such that 
\begin{equation} \label{e.overdet-intro-forms}
	\begin{dcases}
		d \ast d \bE = \lambda \ast \bE,  &\mbox{ in } D\\
		d(\ast \bE)  = 0    &\mbox{ in }   D\\
		i^\ast \bE  = 0 &\mbox{ on } \partial D\\
	i\sqrt{\lambda}	\bH = \ast \,d \bE= 0,&\mbox{ on } \partial D\,, 
	\end{dcases}
\end{equation}
where~$\ast$ denotes the Hodge star operator,~$i : \partial \Omega\to \overline{\Omega}$ is the inclusion map and~$i^\ast$ is the pullback to the boundary. 

Our results in Section~\ref{sec:electro} show that when~$D$ is a ball, then it is possible to write down a solution to this overdetermined problem in terms of spherical harmonics. We end this introduction with

\smallskip 

\textbf{Open question:}
	Does there exist a Lipschitz domain~$D \subset \R^3$ and~$\lambda > 0$ and~$\bE \neq 0$ satisfying~\eqref{e.overdet-intro-forms}  on~$D$, then must~$D$ be a ball?

\subsection{A brief overview of related literature:} There are several papers in the photonics literature that explore devices that one can fabricate using ENZ materials. In what follows we describe some related mathematical work.
\begin{itemize}
	\item The paper~\cite{maier} explores the question of \textit{homogenization} of Maxwell's equations in order to obtain effective tensors that display the ``epsilon near zero" effect. Let us remind the reader that in contrast~\cite{maier}, our paper concerns device \textit{design}-- specifically, designing a resonator that harnesses the effects of an epsilon-near-zero medium. Some related work on homogenization for Maxwell's equations with sign-changing coefficients appears in~\cite{bunoiu2021homogenization}.
	\item More recently, the authors of~\cite{nguyen} study the limiting absorption principle and well-posedness for the time-harmonic Maxwell equations with anisotropic sign-changing coefficients. In our setting, since~$\delta \in \CC$ and we make no assumptions about the sign of the real and imaginary parts of~$\delta,$ our setting, like~\cite{nguyen} has sign-changing coefficients. Our setting, however, departs from that paper since the permittivity~$\e$ degenerates to zero in part of the domain. Some related work in the isotropic setting appears in~\cite{cossetti2021limiting}.
	\item Our previous papers~\cite{KV1,KV2} consider other applications of epsilon-near-zero devices, and our present paper builds on the tools built in those papers-- but take them substantially further in its current application to the setting of a system of equations rather than a single complex scalar equation. The paper~\cite{KV1} dealt with the phenomenon referred to as \textit{photonic doping}; in mathematical terms, we deal there with perturbation theory in the context of a non-self-adjoint \textit{scattering} problem, rather than one of resonance. The paper~\cite{KV2}, on the other hand, is closer to our current paper, since it, too, deals with a spectral perturbation problem, albeit in a simpler scalar setting. 
	\item The recent papers~\cite{Chen-Lipton-ARMA2013,Chen-Lipton-MMS2013,Fortes-Lipton-Shipman-2011,Fortes-Lipton-Shipman-2010} deal with a high-contrast setting similar to ours, in the context of scalar Helmholtz equations; the motivation of these papers lies in photonic band gaps, and is therefore different from our goals. On the other hand, these papers, similar to ours, fall in the realm of spectral perturbation theory, and proceed to prove complex analyticity of the quantities involved by developing their Taylor series and carefully estimating the size of the Taylor coefficients. A particular difference between our approach versus the approach in the aforementioned papers is that we make no use of layer-potential arguments.
\end{itemize}
\section{A warm-up scalar problem}
\label{s.scalar}
In this section we consider a scalar version of the eigenvalue problem in~\eqref{e.strongform}, and analyze it in complete detail. To this goal, let~$\e_\de$ be as in~\eqref{e.ENZmedium} and consider, in place of the vectorial eigenvalue problem~\eqref{e.strongform}, the scalar Neumann eigenvalue problem 
\begin{equation} \label{e.PDEscalar}
	\begin{dcases}
		-\Delta u_\delta &= \lambda_\delta \e_\delta u_\delta\,  \quad \mbox{ in } \Omega\\
		\partial_{\nu_\Omega} u_\de &= 0, \qquad \quad \mbox{ on }\partial \Omega\,. 
	\end{dcases}
\end{equation}
Recall, we think of~$\de \in \CC$ as being small (i.e.,~$|\de|\ll 1$), and so, we will study the~$\de = 0$ version of~\eqref{e.PDEscalar} first. Even though our Maxwell eigenproblem is in the physical dimension~$d=3$ (or, one can formulate a corresponding problem for differential forms in any dimension~$d$), the scalar problem~\eqref{e.PDEscalar} makes sense for~$\Omega,D \subset \Rd$ for general dimension~$d \geqslant 2,$ and our presentation will allow for this slight generality.  
The limiting problem corresponding to~$\de = 0$ is
\begin{equation} \label{e.PDEscalar-0}
	\begin{dcases}
		-\Delta u &= \lambda \indc_D u\,  \quad \mbox{ in } \Omega\\
		\partial_{\nu_\Omega} u &= 0, \qquad \quad \mbox{ on }\partial \Omega\,. 
	\end{dcases}
\end{equation}
\subsection{Functional analytic preliminaries} 
\label{ss:FAprelim}
From the viewpoint of spectral theory, the problems~\eqref{e.PDEscalar} and~\eqref{e.PDEscalar-0} both have different, equivalent, but more convenient reformulations. To explain this reformulation we must introduce some notation. Let~$\{\mu_k:k\in\NN_0\}$ denote the set of Neumann eigenfunctions of~$-\Delta$ in the domain~$\Omega$ with associated eigenfunctions~$\{\chi_k: k \in \NN_0\}$ normalized so that they form an orthonormal basis for~$L^2(\Omega).$  
\begin{equation}
	\begin{dcases}
		-\Delta \chi_k = \mu_k \chi_k \quad \mbox{ in }\Omega\\
		\partial_{\nu_\Omega} \chi_k = 0, \quad \mbox{ on } \partial \Omega\,.
	\end{dcases}
\end{equation}
Here~$\mu_0 = 0$ with~$\chi_0 \equiv \frac1{\sqrt{|\Omega|}},$ and by the assumed connectedness of~$\Omega, $ all other eigenfunctions are nonconstant and have mean-zero. 
Here, and in the rest of the paper, we define 
\begin{equation*}
	\mathring{L}^2(\Omega) := \biggl\{ f \in L^2(\Omega): \int_\Omega f = 0\biggr\}\,. 
\end{equation*}
In light of the connectedness of the domain~$\Omega$, we have the orthogonal decomposition 
\begin{equation*}
	L^2(\Omega) = \mathring{L}^2(\Omega) \oplus \R\,.
\end{equation*}
For any function~$f \in \mathring{L}^2(\Omega),$ 
let us define~$(-\Delta)^{-\sfrac12} f$ spectrally, via 
\begin{equation} \label{e.L-1/2}
	(-\Delta)^{-\sfrac12} f :=  \sum_{k \in \NN} \frac1{\sqrt{\mu_k}} \biggl( \int_\Omega f(y) \chi_k(y)\,dy\biggr)\chi_k\,. 
\end{equation}
Since~$f \in L^2(\Omega),$ this infinite series converges in~$L^2(\Omega),$ as~$\mu_k$ is increasing in~$k,$ and~$\mu_k \to \infty$ as~$k\to \infty.$

\bigskip 
Conversely, let us consider the spectrally defined Sobolev space
\begin{equation*}
	\mathring{H}^{1}(\Omega) := \biggl\{ f \in L^2(\Omega) : \sum_{k=1}^\infty \mu_k \biggl(\int_\Omega f (y) \chi_k(y)\,dy \biggr)^2 < \infty\biggr\}\,;
\end{equation*}
it is a well-known fact that if~$\Omega$ is a bounded, Lipschitz, connected domain, then this coincides with the ``more basic'' definition of the Sobolev space~$\mathring{H}^1(\Omega)$ as mean zero~$L^2$ functions whose distributional gradients admit representation as integration against an~$L^2$ curl-free vector field: 
\begin{equation*}
	\mathring{H}^1(\Omega) \equiv \biggl\{ f \in \mathring{L}^2(\Omega) : \nabla f \in L^2(\Omega)\biggr\}\,,
\end{equation*}
Now, for any~$f \in \mathring{H}^{1}(\Omega)$ 
\begin{equation*}
	(-\Delta)^{\sfrac12} f := \sum_{k=1}^\infty \sqrt{\mu_k} \biggl( \int_\Omega f(y) \chi_k(y)\,dy \biggr)\chi_k\,.
\end{equation*}

We collect several observations in the next lemma. 
\begin{lemma}
	\label{l.observations}
	Let~$\Omega$ be a connected, bounded, Lipschitz domain. Then, 
	\begin{enumerate}
		\item the eigenfunctions~$\{\chi_k\}_{k\in\NN_0}$ with~$\chi_0 \equiv 1$ have the property that the set~$\bigl\{ \frac1{\sqrt{\mu_k}} \chi_k\bigr\}$ are orthonormal in the inner product 
		\begin{equation*}
			\langle u,v \rangle_{\mathring{H}^1(\Omega)} := \int_\Omega \nabla u \cdot \nabla v\,dx\,. 
		\end{equation*}
		\item The operator~$(-\Delta)^{-\sfrac12} :\mathring{L}^2(\Omega) \to \mathring{H}^1(\Omega)$ is an isometric isomorphism so that 
		\begin{equation} \label{e.isom1}
			\| (-\Delta)^{-\sfrac12} f \|_{\mathring{H}^1(\Omega)} = \|f\|_{\mathring{L}^2(\Omega)}\,. 
		\end{equation}
		Similarly, the operator~$(-\Delta)^{\sfrac12}:\mathring{H}^1(\Omega) \to \mathring{L}^2(\Omega)$ is an isometric isomorphism with
		\begin{equation} \label{e.isom2}
			\|(-\Delta)^{\sfrac12} f\|_{\mathring{L}^2(\Omega)} = \| f\|_{\mathring{H}^1(\Omega)}\,. 
		\end{equation}
	\end{enumerate}
\end{lemma}
\begin{proof}
	For the first item, testing the weak formulation of the equation satisfied by~$\chi_j,$ with the test function~$\chi_k,$ we obtain
	\begin{equation*}
		\int_\Omega \nabla \chi_j \cdot \nabla \chi_k \,dx = \mu_j \int_\Omega \chi_j \chi_k\,dx = \mu_j \delta_{j,k}\,.
	\end{equation*}
	For the second item, let~$f \in \mathring{L}^2(\Omega),$ then we can write~$f = \sum_{j=1}^\infty c_j \chi_j$ where the series is convergent in~$\mathring{L}^2(\Omega).$ Then, it is clear that 
	\begin{equation*}
		\|f\|_{\mathring{L}^2(\Omega)}^2 = \sum_{j=1}^\infty c_j^2\,.
	\end{equation*}
	Then, by the definition, 
	\begin{equation*}
		(-\Delta)^{-\sfrac12} f = \sum_{j=1}^\infty \frac{1}{\sqrt{\mu_j}} c_j \chi_j\,.
	\end{equation*}
	Then, by item 1 of the lemma, it follows that 
	\begin{equation*}
		\|f\|_{\mathring{H}^1(\Omega)}^2 = \sum_{j=1}^\infty \mu_j \frac1{\mu_j} c_j^2 = \sum_{j=1}^\infty c_j^2 = \|f\|_{\mathring{L}^2(\Omega)}^2\,. 
	\end{equation*}
	This is~\eqref{e.isom1}. The proof of the assertion~\eqref{e.isom2} is identical.
\end{proof}
\begin{remark}
	\label{r.weakvsstrong} If~$f \in \mathring{L}^2(\Omega),$ then it is well-known that there is a unique~$u \in \mathring{H}^1(\Omega)$ that is a weak solution to 
	\begin{equation*}
		\begin{dcases}
			-\Delta u = f \quad \mbox{ in } \Omega\\
			\partial_\nu u = 0 \quad \mbox{ on } \partial \Omega\,,
		\end{dcases}
	\end{equation*}
	that is, 
	\begin{equation*}
		\int_\Omega \langle  \nabla u , \nabla \phi \rangle \,dx= \int_\Omega f \phi \,dx \quad \mbox{ for all } \phi \in \mathring{H}^1(\Omega)\,.
	\end{equation*}
	It is also then immediate that if~$f = \sum_{j=1}^\infty c_j \chi_j,$ then~$u = \sum_{j=1}^\infty \frac{c_j}{\mu_j} \chi_j\,.$ This infinite series is understood as a convergent series in~$\mathring{H}^1(\Omega)$, since by the first item in Lemma~\ref{l.observations} we compute 
	\begin{multline*}
		\|u\|_{\mathring{H}^1(\Omega)}^2 =  \int_\Omega  \nabla u \cdot \nabla u \,dx = \sum_{j=1}^\infty \frac{c_j^2}{\mu_j^2} \int_\Omega |\nabla \chi_j|^2\,dx = \sum_{j=1}^\infty \frac{c_j^2}{\mu_j^2}\mu_j \\
		= \sum_{j=1}^\infty \frac{c_j^2}{\mu_j} \leqslant \frac1{\mu_1} \sum_{j=1}^\infty c_j^2 = \frac1{\mu_1} \|f\|_{\mathring{L}^2(\Omega)}^2\,. 
	\end{multline*}
\end{remark}
Finally, we have
\begin{lemma} \label{l.equivalence-2}
	If~$f \in \mathring{L}^2(\Omega)$ and~$w \in \mathring{H}^1(\Omega)$ satisfy 
	\begin{equation} \label{e.equality-1}
		(-\Delta)^{-\sfrac12} f = (-\Delta)^{\sfrac12} w\,,
	\end{equation}
	then, in fact, 
	\begin{equation} \label{e.weak solution}
		\begin{dcases}
			-\Delta w = f \mbox{ in } \Omega\\
			\partial_\nu w = 0 \mbox{ on } \partial \Omega\,,
		\end{dcases}
	\end{equation}
	in the weak sense. 
\end{lemma}
\begin{proof}
	Writing~$f \in \mathring{L}^2(\Omega)$ in the orthonormal basis~$\{\chi_j\}_{j=1}^\infty$, i.e.,if~$f = \sum_{j=1}^\infty c_j \chi_j,$ for some~$\{c_j\}_{j=1}^\infty $ then~$(-\Delta)^{-\sfrac12} f = \sum_{j=1}^\infty \frac1{\sqrt{\mu_j}} c_j \chi_j = \sum_{j=1}^\infty \sqrt{\mu_j}\frac1{\mu_j} c_j \chi_j = (-\Delta)^{\sfrac12} w$ with~$w = \sum_{j=1}^\infty \frac1{\mu_j} c_j\chi_j.$ But by Remark~\ref{r.weakvsstrong} then,~\eqref{e.weak solution} holds. This entails, in particular, the homogeneous Neumann boundary condition~$\partial_\nu w = 0.$ 
\end{proof}

\subsection{Analysis of the limiting scalar problem} \label{ss.scalarlimit}

Having set up this functional analytic framework, we are ready to study the generalized eigenvalue problem~\eqref{e.PDEscalar-0} through an indirect argument. For this, we introduce the operator~$K_0$ on~$\mathring{L}^2(\Omega)$ via 
\begin{equation}
	\label{e.K_0def}
	K_0 : \mathring{L}(\Omega) \to \mathring{L} (\Omega) \quad K_0 v := (-\Delta)^{-\sfrac12} \indc_D (-\Delta)^{-\sfrac12} v\,.
\end{equation}
We may extend this definition to~$L^2(\Omega)$ as follows: for any~$f \in L^2(\Omega)$, we let~$\P_0 (f) := f - \fint_\Omega f$ denote the projection to the orthogonal complement of the constants, namely, mean-zero functions. Then, we extend~$K_0$ by setting 
\begin{equation} \label{e.projectiondefK_0}
	K_0 := \P_0^T (-\Delta)^{-\sfrac12} \P_0 \indc_D \P_0^T (-\Delta)^{-\sfrac12}\P_0\,.
\end{equation}
Below in Lemma~\ref{l.K0good} we show that~$K_0$ is a compact, self-adjoint, nonnegative operator on~$L^2(\Omega).$ By the spectral theorem for compact, self-adjoint operators we infer that the the spectrum of~$K_0$ consists of a discrete sequence of positive eigenvalues that can be arranged in non-increasing order, repeated per (finite) multiplicity that accumulate at~$0,$ and in addition the associated eigenfunctions can be normalized to form an orthonormal Hilbert basis for the orthogonal compliment of~$\ker (K_0)$ in~$L^2(\Omega)$. 


Let now~$v$ be one such eigenfunction, with associated eigenvalue~$\rho > 0.$ Defining
\begin{equation*}
	\lambda := \rho^{-1}, \quad \mbox{and } \tilde{u} := (-\Delta)^{-\sfrac12} v\,,
\end{equation*} 
the preceding discussion implies that~$\tilde{u}$ satisfies 
\begin{equation*}
	\Delta \tilde{u} - \lambda \indc_D\tilde{u}  = \mbox{ constant in }  \Omega\,, \quad \partial_{\nu_\Omega} \tilde{u} = 0 \quad \mbox{ on }\partial \Omega\,.
\end{equation*}
Integrating both sides and utilizing the boundary condition, we find that the required constant is given by~$\lambda \frac{1}{|\Omega|} \int_D \tilde{u},$ so that setting~$u := \tilde{u} - \frac1{|\Omega|}\int_D \tilde{u}$ yields an eigenfunction of the generalized eigenvalue problem in~\eqref{e.PDEscalar-0}. 

It is also clear that this argument can be reversed, and therefore this procedure yields every nontrivial eigenpair to the generalized eigenvalue problem~\eqref{e.PDEscalar-0}.

\smallskip
We turn to detailing and formalizing the foregoing discussion. The first lemma contains properties of the operator~$K_0.$ 
\begin{lemma}
	\label{l.K0good}
	The operator~$K_0$ defined in~\eqref{e.K_0def} is a compact, self-adjoint, nonnegative operator on~$L^2(\Omega)$ to itself. 
\end{lemma}
\begin{proof}
	Self-adjointness of~$K_0.$ In order to study~$K_0,$ we may always suppose it acts on~$f \in \mathring{L}^2(\Omega)$ since~$K_0$ begins by projection to mean-zero functions. 
	
	\smallskip The compactness of~$K_0$ follows from the compactness of~$(-\Delta)^{-\sfrac12}$ and the boundedness of multiplication by~$\indc_D.$ The boundedness of the latter is easy to verify. The compactness of~$(-\Delta)^{-\sfrac12}$ follows as we can approximate it in the operator norm by a sequence of finite rank operators: for any~$M\in \NN,$ let us set 
	\begin{equation*}
		P_M (-\Delta)^{-\sfrac12} : \mathring{L}^2(\Omega) \to \mathring{L}^2(\Omega) := \sum_{k=1}^{M} \frac{1}{\sqrt{\mu_k}} \biggl( \int_{\Omega} f(y)\chi_k(y)\,dy\biggr)\chi_k\,.
	\end{equation*} 
	It is clear that~$P_M (-\Delta)^{-\sfrac12}$ is a finite rank operator. Also note that by definition
	\begin{equation*}
		\|P_M (-\Delta)^{-\sfrac12} f  - (-\Delta)^{-\sfrac12} f \|_{L^2(\Omega)} \leqslant \frac1{\sqrt{\mu_{M+1}}} \|f\|_{L^2(\Omega)}\,. 
	\end{equation*}
	But~$\mu_M\to \infty$ as~$M \to \infty,$ so the right-hand side of the preceding equation goes to zero as we send~$M \to \infty.$ Therefore~$(-\Delta)^{-\sfrac12} f$ is compact, and hence,~$K_0$ is, too. 
	
	\smallskip 
	To see the assertion about the nonnegativity of~$K_0$, let us compute that
	\begin{equation*}
		\int_\Omega (K_0 v ) v = \int_\Omega \indc_D|(-\Delta)^{-\sfrac12} v |^2 \,dx = \int_D|(-\Delta)^{-\sfrac12} v|^2\geqslant 0\,.
	\end{equation*}
\end{proof}
The spectral theorem for compact self-adjoint operators then asserts that the spectrum of~$K_0$ is discrete and is comprised of a sequence~$\{\rho_k: k \in \NN\}$ with~$\rho_{k+1} \leqslant \rho_k,$ converging to zero as~$k \to \infty;$ in addition each~$\rho_k >0 $ has finite multiplicity, and the associated eigenfunctions form an orthonormal Hilbert basis for~$\ker(K_0)^\perp$. The next lemma characterizes the kernel of~$K_0.$  
\begin{lemma}
	The kernel of the operator~$K_0$ is given by 
	\begin{equation*}
		\ker (K_0) :=\Bigl \{v \in L^2(\Omega) :  v = (-\Delta)^{\sfrac12} g \mbox{ for some } g \in L^2(\Omega) \mbox{ with } g = 0 \mbox{ a.e. on } D\Bigr\}\,. 
	\end{equation*}
\end{lemma}
\begin{proof}
	Let~$v \in \ker (K_0),$ then as computed in the previous lemma we find
	\begin{equation*}
		0 = \int_\Omega (K_0 v)v \,dx = \int_D|(-\Delta)^{-\sfrac12} v|^2\,dx \,,
	\end{equation*}
	so that we find 
	\begin{equation*}
		(-\Delta)^{-\sfrac12} v = 0, \quad \mbox{ a.e. in } D\,.
	\end{equation*}
	Set~$(-\Delta)^{-\sfrac12} v = g,$ so that~$g = 0$ a.e. on~$D.$ Then, for each~$k \in \NN$
	\begin{equation*}
		\int_\Omega g (y)\chi_k(y) \,dy = \int_{\ENZ} g(y) \chi_k(y)\,dy = \frac1{\sqrt{\mu_k}} \int_\Omega v (y) \chi_k(y)\,dy\,, 
	\end{equation*}
	defines the projection of~$g$ up to its average on~$\Omega,$ which~$(-\Delta)^{\sfrac12}$ does not see. 
\end{proof}
Putting together the ingredients above, we have thus proved:
\begin{theorem}
	\label{t.limscalar}
	There exist a countable number of generalized eigenvalues~$\{\lambda_k\}_{k=1}^\infty$ with~$0 < \lambda_k \to \infty$ as~$k \to \infty$, and associated generalized eigenfunctions~$\{\phi_k\}_{k=1}^\infty$ satisfying 
	\begin{equation*}
		\begin{dcases}
			-\Delta \phi_k = \la_k \indc_D \phi_k, \quad &\mbox{ in }\Omega\\
			\partial_{\nu_\Omega} \phi_k = 0 \quad & \mbox{ on }\partial \Omega\,. 
		\end{dcases}
	\end{equation*}
	These eigenfunctions satisfy~$\int_D \phi_k = 0$ and are orthogonal according to 
	\begin{equation*}
		\int_D\phi_k \phi_l \,dx = 0, \quad k \neq l\,. 
	\end{equation*}
\end{theorem}
\subsection{Perturbation theory and analysis of~\eqref{e.PDEscalar}}
We turn our attention next to an analysis of~\eqref{e.PDEscalar}, and we need some functional analytic framework. Consider the bounded linear functional~$\ell_\de$ on~$\mathring{L}^2(\Omega)$ defined as follows: for any~$v \in \mathring{L}^2(\Omega)$, ~$\ell_\de(v)$ denotes the unique constant which is such that 
\begin{equation*}
	\int_\Omega \e_\de \Bigl( (-\Delta)^{-\sfrac12} v + \ell_\de(v) \Bigr)\,dx = 0\,. 
\end{equation*}
Explicitly, 
\begin{equation*}
	\ell_\de (v) := - \frac{ \int_\Omega \e_\de (-\Delta )^{-\sfrac12} v \,dx }{\int_\Omega \e_\de \,dx} = \frac{\int_D (-\Delta)^{-\sfrac12} v\,dx + \de \int_{\ENZ} (-\Delta)^{-\sfrac12} v\,dx }{ |D| + \de |\ENZ|}\,. 
\end{equation*}
The following lemma is immediate: 
\begin{lemma} \label{l.analdepofldel}
	The family of linear functionals~$\{\ell_\de\}_{\de}$, for~$|\de| < \frac{|D|}{|\ENZ|},$ are continuous linear functionals on~$\mathring{L}^2(\Omega)$ that depend analytically on~$\de$ for~$|\de| < \frac{|D|}{|\ENZ|}. $
\end{lemma}
Towards analyzing the system~\eqref{e.PDEscalar}, we define the operator~$K_\de$ on~$\mathring{L}^2(\Omega)$, defined by 
\begin{equation} \label{e.Kdelta}
	K_\de := (-\Delta)^{-\sfrac12} \e_\de \Bigl( (-\Delta)^{-\sfrac12} + \ell_\de \Bigr)   \,,
\end{equation}
In the next lemma, we record the basic properties of~$K_\de$ that are crucial subsequently. 
\begin{lemma}
	\label{l.propertiesofKdelta}
	The family of operators~$K_\de: \mathring{L}^2(\Omega)\to \mathring{L}^2(\Omega)$ is compact for every~$\de \in \CC$ with~$|\de| < \frac{|D|}{|\ENZ|}$, that depend analytically on~$\de$ in this neighborhood of~$\de = 0.$ In addition, for~$\de \in \RR$ in this neighborhood, the operator~$K_\de$ is also self-adjoint.
\end{lemma}
\begin{proof}
	The proof of compactness proceeds identically as in Lemma~\ref{l.K0good} and is omitted. Analytic dependence of~$K_\de$ on~$\de$ is also immediate, and therefore we proceed to check that~$K_\de$ is self-adjoint when~$\de \in \RR$ with~$|\de| < \frac{|D|}{|\ENZ|}$. To see this, for any~$v,w \in \mathring{L}^2(\Omega)$ we compute   
	\begin{multline}
		\langle K_\de v, w \rangle_{\mathring{L}^2(\Omega)} = \int_\Omega \Bigl( (-\Delta)^{-\sfrac12} \e_\de \Bigl( (-\Delta)^{-\sfrac12} v + \ell_\de(v) \Bigr) \Bigr) w \\
		= \int_\Omega \e_\de \Bigl( (-\Delta)^{-\sfrac12} v + \ell_\de(v) \Bigr)  (-\Delta)^{-\sfrac12} w \\
		= \int_\Omega \e_\de \Bigl( (-\Delta)^{-\sfrac12} v + \ell_\de(v) \Bigr)  \Bigl( (-\Delta)^{-\sfrac12} w + \ell_\de (w) \Bigr)\,dx \,. 
	\end{multline}
	In the last line, we used the fact that~$\ell_\de(w)$ is constant, and that~$\ell_\de(v)$ is chosen to make~$\e_\de \Bigl( (-\Delta)^{-\sfrac12} v  + \ell_\de(v) \Bigr)$ orthogonal to constants. The last line of the display is symmetric in~$v,w$ and therefore we can switch the roles of~$v,w$ to obtain self-adjointness in this case. 
\end{proof}
It follows by the spectral theorem for compact operators that for each~$\de$ with~$|\de|<\frac{|D|}{|\ENZ|},$ the operator~$K_\de$ admits a sequence of eigenvalues~$\rho_\de,$ that are real when~$\de$ is real (by self-adjointness). In addition, the associated eigenfunctions form a Hilbert basis for the orthogonal compliment of~$\ker (K_\de).$ 

The next proposition is crucial in connecting the spectrum of the operator~$K_\de$ to eigenpairs satisfying~\eqref{e.PDEscalar}. 
\begin{proposition}
	\label{p.equivalence of spectral problems}
	The following statements are equivalent: 
	\begin{enumerate}
		\item The number~$\la_\de \neq 0$ is an eigenvalue of~\eqref{e.PDEscalar}, for which there exists an associated function~$u_\de \in H^1(\Omega),$ where the pair~$(u_\de,\la_\de)$ solves the generalized eigenvalue problem~\eqref{e.PDEscalar}:
		\begin{equation} \label{e.v1}
			\begin{dcases}
				-\Delta u_\de = \la_\de \e_\de u_\de \quad \mbox{ in } \Omega\\
				\partial_{\nu} u_\de = 0, \quad \mbox{ on } \partial \Omega\,. 
			\end{dcases}
		\end{equation}
		\item The number~$\la_\de^{-1} \neq 0$ is an eigenvalue of the compact operator~$K_\de$, i.e., there exists an eigenfunction~$v_\de \in \mathring{L}^2(\Omega)$ which solves 
		\begin{equation} \label{e.v2}
			K_\de v_\de = \la_\de^{-1} v_\de\,.
		\end{equation}
	\end{enumerate}
	Furthermore, if one of these (and therefore, both of them) holds, then 
	\begin{equation} \label{e.udelvdel}
		u_\de = (-\Delta)^{-\sfrac12} v_\de + \ell_\de(v_\de)\,. 
	\end{equation}
\end{proposition}
\begin{proof}
	\underline{Proof of (1) ~$implies$ (2):} Given~$u_\de \in H^1(\Omega)$ solving~\eqref{e.v1}, then there exists a unique~$v_\de \in \mathring{L}^2(\Omega)$ such that 
	\begin{equation*}
		u_\de = (-\Delta)^{-\sfrac12} v_\de + \ell_\de (v_\de)\,.
	\end{equation*}
	Thanks to~\eqref{e.v1} we have 
	\begin{equation*}
		-\Delta u_\de =   -\Delta \Bigl( (-\Delta)^{-\sfrac12} v_\de + \ell_\de(v_\de)\Bigr) = \la_\de \e_\de  \Bigl( (-\Delta)^{-\sfrac12} v_\de + \ell_\de (v_\de) \Bigr)\,. 
	\end{equation*}
	But~$\ell_\de(v_\de)$ is a constant, and~$(-\Delta)^{\sfrac12} v_\de \in \mathring{L}^2(\Omega).$ Applying~$(-\Delta)^{-\sfrac12}$ to both sides yields 
	\begin{equation*}
		v_\de  = \la_\de K_\de(v_\de)\,,
	\end{equation*}
	which is~\eqref{e.v2}. 
	
	\medskip 
	\underline{Proof of (2) implies (1):} Given~$v_\de \in \mathring{L}^2(\Omega),$ we let~$u_\de$ be given by~\eqref{e.udelvdel}. Plugging this in~\eqref{e.v2} we obtain 
	\begin{multline} \label{e.Kdeltaegfct}
		K_\de v_\de = (-\Delta)^{-\sfrac12} \e_\de \Bigl( (-\Delta)^{-\sfrac12} v_\de + \ell_\de(v_\de) \Bigr)  = (-\Delta)^{-\sfrac12} \e_\de u_\de = \la_\de^{-1} v_\de \\
		= \la_\de^{-1} (-\Delta)^{\sfrac12} \Bigl( u_\de - \ell_\de(v_\de)\Bigr)\,. 
	\end{multline}
	Applying Lemma~\ref{l.equivalence-2} with~$f = \e_\de u_\de$ and~$w = \la_\de^{-1} (u_\de - \ell_\de(v_\de)),$ we obtain 
	\begin{equation*}
		\begin{dcases}
			-\Delta w = -\la_{\de}^{-1} \Delta u_\de = \e_\de u_\de \quad \mbox{ in } \Omega\\
			\partial_{\nu} w = \la_{\de}^{-1} \partial_\nu u_\de = 0, \quad \mbox{ on } \partial \Omega\,. 
		\end{dcases}
	\end{equation*}
	Rearranging this yields~\eqref{e.v1}. 
	
	This completes the equivalence of (1) and (2), and also, the relationship between~$u_\de$ and~$v_\de.$ 
\end{proof}

Turning to the dependence of the spectrum of~$K_\de$ on~$\de,$ we see that it is a family of self-adjoint, compact operators on~$\mathring{L}^2(\Omega)$ that depend analytically on~$\de$ about~$\de = 0.$ Rellich's theorem, recalled in~Theorem \ref{t.rellich}, immediately applies. To state the conclusions that we can draw from Rellich's theorem, let us introduce some notation. First, let~$\la_{0,k}^{-1}$ denote the~$k$th eigenvalue of~$K_0$ introduced in~\eqref{e.K_0def} and let~$\la_{\de,k}^{-1}$ denote the~$k$th eigenvalue of~$K_\de$ introduced in~\eqref{e.Kdelta}. Let~$\gamma(\la_{0,k})$ denote the spectral gap about~$\la_{0,k}:$
\begin{equation*}
	\ga(\la_{0,k}) := \min_{j \in \NN} \{|\la_{0,j} - \la_{0,k}| : \la_{0,j} \neq \la_{0,k}\}>0\,. 
\end{equation*}
The main theorem of this section is:
\begin{theorem}
	\label{t.scalar}
	For each~$k \in \NN$, let~$\lambda_{0,k}$ have (finite) multiplicity~$h,$ and let~$0 < \gamma < \gamma(\la_{0,k})$ be fixed, so that~$K_0$ has no eigenvalues in closed neighborhood~$[\gamma(\la_{0,k}) - \gamma, \gamma(\la_{0,k}) + \gamma]$ other than~$\la_{0,k}$, and suppose, without loss of generality, that~$\la_{0,k-1} < \la_{0,k} = \la_{0,k+1} = \cdots = \la_{0,k+h-1} < \la_{0,k+h}.$ Then there exists a positive number~$\delta_0 = \delta_0(k) \leqslant \frac{|D|}{|\ENZ|},$ and~$h$ families of complex numbers 
	\begin{equation*}
		\la_{\de,k}, \cdots, \la_{\de,k+h-1}\,,
	\end{equation*}
	that have convergent power series in~$\de$ in a neighborhood~$|\de|< \de_0(k)$; additionally, there exist an associated family of functions~$\{u_{\de,k},\cdots, u_{\de,k+h-1}\} \subset H^1(\Omega)$ also given by a convergent power series in the Hilbert space~$H^1(\Omega)$ that satisfy the following properties: 
	For each~$i = 0,\cdots, h-1,$ we have 
	\begin{equation} \label{e.udelk+i}
		\begin{dcases}
			- \Delta u_{\de,k+i} &= \la_{\de,k+i} \e_\de u_{\de,k+i}, \quad \mbox{ in } \Omega\\
			\partial_{\nu} u_{\de,k+i} &= 0, \quad \quad \mbox{ on } \partial \Omega\,.
		\end{dcases}
	\end{equation}
	Furthermore, for each~$i \in \{0,\cdots, h-1\}$, we have that~$\lim_{\de \to 0}\la_{\de,k+i} = \la_{0,k+i},$ and in addition, for all~$\de, |\de|<\de_0(k),$ and~$i,j = 0,\cdots,h-1,$ the following orthogonality relation holds:
	\begin{equation*}
		\int_{\Omega}\e_\de u_{\de,k+i} u_{\de,k+j} = \de_{ij} = \begin{dcases}
			1 & \mbox{ if } i = j\\
			0 & \mbox{ otherwise } 
		\end{dcases}\,. 
	\end{equation*}
	Finally, for~$|\de| \ll 1,$ the interval~$[\la_{0,k} - \frac{\gamma}2, \la_{0,k} + \frac{\gamma}2]$ includes no other eigenvalues for the generalized eigenvalue problem~\eqref{e.PDEscalar} except for the~$h$ branches~$\la_{\de,k}, \cdots, \la_{\de,k+h-1}.$ 
\end{theorem}
\begin{proof}
	By Lemmas~\ref{l.analdepofldel} and~\ref{l.propertiesofKdelta}, the family of operators~$\de \mapsto K_\de$ is a family of compact operators on~$\mathring{L}^2(\Omega),$ which are self-adjoint for~$|\de| \ll 1$ when~$\de \in \RR.$ Therefore Rellich theorem applies, and immediately yields~$h$ analytic branches~$\la_{\de,k}^{-1}, \cdots, \la_{\de,k+h-1}^{-1}$ of eigenvalues and associated eigenfunctions~$\{v_{\de,k}, \cdots , v_{\de,k+h-1}\}$ of the operator~$K_\de$ that all admit convergent power series representations in~$\de$ for~$|\de| < \de_0(k)$ sufficiently small. In addition, these satisfy~$\lim_{\de \to 0} \la_{\de,k+i} = \la_{0,k+i}$ for each~$i=0,\cdots,h-1.$ Since~$\la_{0,k+i} > 0$ for each~$i,$ it follows that the maps~$\de \mapsto \la_{\de,k}$ are also analytic. Finally, by Proposition~\ref{p.equivalence of spectral problems}, we can construct eigenfunctions~$u_{\de,k+i}$ associated to~\eqref{e.PDEscalar} with eigenvalues~$\la_{\de,k+i}$ using the relation
	\begin{equation*}
		u_{\de,k+i} := (-\Delta)^{-\sfrac12} v_{\de,k+i} + \ell_\de(v_{\de,k+i})\,, \quad i = 0,\cdots, h-1\,.
	\end{equation*}
	Since by Lemma~\ref{l.analdepofldel} the dependence of~$\ell_\de$ is analytic on~$\de,$ and by Rellich's theorem the dependence of each~$v_{\de,k+i}$ is analytic on~$\de,$ the same is true of each~$u_{\de,k+i}.$ 
	Finally, we must check the orthogonality condition. If~$\lambda_{\de,k+i} \neq \lambda_{\de,k+j}$ then this follows from multiplying the equation for~$u_{\de,k+i}$ by~$u_{\de,k+j}$ and vice-versa, integrating, and subtracting the results:
	\begin{equation*}
		0 = \int_{\Omega } \nabla u_{\de,k+i} \cdot \nabla u_{\de,k+j} \,dx = \bigl( \la_{\de,k+i} - \la_{\de,k+j}\bigr) \int_\Omega \e_\de u_{\de,k+i}u_{\de,k+j}\,dx\,. 
	\end{equation*}
	If it happens that~$\la_{\de,k+i} = \la_{\de,k+j}$ also for (some)~$\de \neq 0$ (as it might, rarely), then there is a different argument to prove the same orthogonality: using~\eqref{e.Kdeltaegfct}, the facts that~$\ell_\de(v_{\de,k+i})$ is constant, and that~$\ell_\de(v_{k+j})$ is chosen to make~$\e_\de u_{\de,k+j}$ orthogonal to constants, and finally, the self-adjointness of~$(-\Delta)^{-\sfrac12},$ we compute
	\begin{multline}
		\int_\Omega \e_\de u_{\de,k+i} u_{\de,k+j}\,dx \\
		= \int_\Omega \e_\de \Bigl( (-\Delta)^{-\sfrac12} v_{\de,k+i} + \ell_\de (v_{\de,k+i}) \Bigr) \Bigl( (-\Delta)^{-\sfrac12} v_{\de,k+j} + \ell_\de (v_{\de,k+j}) \Bigr)\,dx \\
		= \int_\Omega \e_\de  (-\Delta)^{-\sfrac12} v_{\de,k+i} \Bigl( (-\Delta)^{-\sfrac12} v_{\de,k+j} + \ell_\de (v_{\de,k+j}) \Bigr)\,dx \\
		= \int_\Omega   (-\Delta)^{-\sfrac12} v_{\de,k+i} \e_\de\Bigl( (-\Delta)^{-\sfrac12} v_{\de,k+j} + \ell_\de (v_{\de,k+j}) \Bigr)\,dx \\
		= \int_\Omega v_{\de,k+i} (-\Delta)^{-\sfrac12} \e_\de\Bigl( (-\Delta)^{-\sfrac12} v_{\de,k+j} + \ell_\de (v_{\de,k+j}) \Bigr)\,dx\\
		= \la_{\de,k+j}^{-1} \int_\Omega v_{\de,k+i} v_{\de,k+j}\,dx = 0\,,
	\end{multline}
	since~$v_{\de,k+i}$ and~$v_{\de,k+j}$ are~$L^2(\Omega)-$orthogonal by virtue of being eigenfunctions of the self-adjoint operator~$K_\de.$ 
\end{proof}

\section{Spectral perturbation theory for the Maxwell eigenproblem} \label{sec:vec}
Next we turn to the analysis of~\eqref{e.strongform}. Its limiting version corresponding to~$\de = 0$ is 
\begin{equation} \label{e.vectorde=0}
	\nabla \times \nabla \times \bE_0 = \la_0 \indc_D\bE_0 \quad  \mbox{ in } \Omega, \qquad \bE_0 \times \nu_\Omega = 0, \quad  \mbox{ on } \partial \Omega\,,
\end{equation}
and so we will study the equation~\eqref{e.strongform} perturbatively about the spectrum of the preceding equation. This section is organized in three sub-sections: in the first subsection we gather properties of a family of Helmholtz decompositions that depend on~$\de$. Then, in subsection we study the limiting problem corresponding to~$\de = 0$ through a related compact operator~$\mathbb{K}_0.$ Finally, in the last sub-section we study the problem for~$\de \neq 0$ by putting the spectrum of~\eqref{e.strongform} in  (analytically varying) bijective correspondence with that of a family of compact operators~$K_\de$ on a suitable Hilbert space, that also depend analytically on~$\de$ for~$|\de| \ll 1$. It is to this family of operators that Rellich's theorem applies. 

\subsection{Some Helmholtz Projections}
\label{s.helmholtz}
Throughout this section we will repeatedly work with a particular Helmholtz projection and so we record basic facts about it, for convenience. 

In light of the vanishing of the tangential boundary condition of the electric field~$\bE_0$ in~\eqref{e.vectorde=0}, the convenient form of Helmholtz projection is as follows: for any vector field~$\bF \in L^2(\Omega)^3$ we consider the decomposition of the form~$\bF = \nabla h + w$ where
\begin{equation*}
	h \in H^1_c(\Omega) := \bigl\{ h \in H^1(\Omega): h = \mbox{ constant on } \partial \Omega\bigr\}\,. 
\end{equation*}
This decomposition is unique, and results from solving a least squares problem: the function~$h$ is given by 
\begin{equation*}
	h = \argmin \biggl\{\int_\Omega |\bF - \nabla q|^2\,dx : q \in H^1_c(\Omega)\biggr\}\,.
\end{equation*}
The minimizer~$h$ is unique by strict convexity, and then~$w := \bF - \nabla h.$ From the criticality condition we know that 
\begin{equation*}
	\mathrm{div}\, w = 0\, \quad \mbox{ in } \Omega, \quad \int_{\partial \Omega}  (\bF - \nabla h) \cdot \nu \,d \sh^2 = 0\,. 
\end{equation*}
Let us note that if~$\bF \in H(\mathrm{curl}, \Omega)$ with~$\bF \times \nu_\Omega = 0$ on~$\partial \Omega,$ then the pieces of this Helmholtz decomposition inherit this property: since~$h$ is constant on~$\partial \Omega$ its tangential gradient is zero on the boundary, and therefore also~$w \times \nu_\Omega = 0$ on~$\partial \Omega.$

Throughout the sequel we use the projection operator~$\P_\Omega :L^2(\Omega)^3 \to L^2(\Omega)^3$ defined by
\begin{equation}
	\label{e.Pomegadef}
	\P_\Omega \bF = w\,,
\end{equation}
with~$w$ as above.

\medskip 
We also are led to consider a family of~$\de$ dependent Helmholtz projections of a fixed divergence free vector field, and must show that these, too, depend analytically on~$\de$. It is a special case of the following more general lemma that we will require, where we consider a family of vector fields~$\mathbf{F}_\de$ that depend analytically on~$\de,$ and show that they have certain~$\de$- dependent Helmholtz projections that depend on~$\de$ analytically. The proof of this lemma, which we believe could be of independent interest, is proven by a strategy similar to~\cite{KV1}. 

\begin{lemma}
	\label{l.analyticity-poisson-highcontrast} 
	Let~$\{\mathbf{F}_{\de}\}_{\de \in \mathbb{C}: |\de|<\de_0} \subset H_0(\mathrm{curl}, \Omega)$ of divergence-free vector fields that depend analytically on~$\de$. Then, there exists a unique solution~$h$ solving
	\begin{equation} \label{e.poisson-statement}
		\begin{dcases}
			\mathrm{div}\, \e_\de\bigl(\mathbf{F}_{\de} + \nabla h \bigr) = 0, \quad & \mbox{ in }\Omega\\
			h= \mathrm{const.} \quad & \mbox{ on }\partial \Omega\\
			\int_{\partial \Omega} \bigl( \bF_\de + \nabla h\bigr) \cdot \nu_\Omega = 0\,. 
		\end{dcases}
	\end{equation}
	Moreover, the function~$h \in H^1_c(\Omega)$ depends analytically on~$\de$ near~$\de=0.$ 
\end{lemma}
\begin{proof}
	We organize the proof in several steps. 
	
	\medskip 
	\textbf{Step 0. } In this step we set up the structure of the proof. 
	Since~$\bF_{\de}$ depends analytically on~$\de$ (taking values in the Hilbert space~$H:=  H_0(\mathrm{curl}, \Omega)$) in a neighborhood of ~$\de=0,$ let us write 
	\begin{equation*}
		\bF_{\de} = \sum_{k=0}^\infty \de^k\bF_{k}\,, 
	\end{equation*}
	for some vector fields~$\bF_{k}$ that do not depend on~$\de,$ and crucially, are \textbf{divergence-free} in~$\Omega.$ 
	This power series converges in the norm of the above Hilbert space~$H$. 
	Let us look for a solution~$h$ to~\eqref{e.poisson-statement} in the form 
	\begin{equation*}
		h = \sum_{k=0}^\infty \de^k h_k\,,
	\end{equation*}
	for some functions~$h_k \in H^1_c(\Omega)$ independent of~$\de$, and construct these by induction on~$k.$  
	
	\medskip 
	\textbf{Step 1.} In this step we introduce certain operators that will be used in the construction, and estimate their operator norms. 
	
	For any distribution (normal flux)~$\sigma$ with~$\sigma \cdot \nu_D\in H^{-\sfrac12}(\partial D)$ which is such that
	\[\int_{\partial D} \sigma \cdot \nu \,d\sh^2 = \langle \sigma \cdot \nu , 1 \rangle_{H^{-\sfrac12} \times H^{\sfrac12} (\partial D)} = 0\,,\]
	we let~$h_\sigma$ denote the unique \textbf{mean-zero} solution to the Neumann problem 
	\begin{equation*}
		\begin{dcases}
			-\Delta h_\sigma = 0, \quad \mbox{ in } D\\
			\partial_{\nu_D} h_\sigma = \sigma\cdot \nu_D\, \quad \mbox{ on } \partial D\,. 
		\end{dcases}
	\end{equation*}
	Let us denote the map~$T_D :\sigma \cdot\nu_D \in H^{-\sfrac12}(\partial D) \mapsto h_\sigma \in H^1(D),$ and observe that~$T_D$ is a bounded linear map by elliptic theory \cite{evans}, and so we find 
	\begin{equation} \label{e.TDbound}
		\|T_D \sigma\|_{H^1(D)} \leqslant C_D \|\sigma\cdot \nu\|_{H^{-\sfrac12}(\partial D)}\,.
	\end{equation}
	Similarly, given~$g \in H^{\sfrac12}(\partial D)$ let~$\mathring{h}_g \in H^1(\ENZ)$ denote the unique solution to
	\begin{equation*}
		\begin{dcases}
			\Delta \mathring{h}_g = 0 \quad \mbox{ in } \ENZ\\
			\mathring{h}_g = g \quad \mbox{ on } \partial D\\
			\mathring{h}_g = 0 \quad \mbox{ on } \partial \Omega\,.
		\end{dcases} 
	\end{equation*}
	Denote by~$T_{\ENZ}: g \in H^{\sfrac12}(\partial D) \mapsto \mathring{h}_g \in H^1(\ENZ),$ which satisfies 
	\begin{equation} \label{e.TENZbound}
		\|\mathring{h}_g\|_{H^1(\ENZ)} \leqslant C_{\ENZ}\|g\|_{H^{\sfrac12}(\partial D)}\,. 
	\end{equation}
	We let~$\gamma_{0,D}$ denote the usual trace operator that associates a function~$h \in H^1(D)$ to its boundary trace~$\gamma_{0,D}(h) \in H^{\sfrac12}(\partial D).$ In particular, if~$h$ is a harmonic function on~$D$, then \footnote{Essentially, one can take this as a definition of the~$H^{\sfrac12}$ norm of a function~$g$ on~$\partial D:$ namely, the Dirichlet energy of the minimizer of~$\int_D |\nabla v|^2\,dx$ among competitors~$v = g$ at~$\partial D.$ If one takes a different definition, then~\eqref{e.bdrytrace} is an equivalence, rather than an exact equality. This distinction will not matter in the sequel. }
	\begin{equation}
		\label{e.bdrytrace}
		\|\gamma_{0,D}(h)\|_{H^{\sfrac12}(\partial D)} = \|\nabla h\|_{L^2(D)}\,. 
	\end{equation}
	Finally, we consider the linear functional given by: for any~$\sigma$ with~$\sigma\cdot \nu_\Omega \in H^{-\sfrac12}(\partial \Omega)$ we let~$\mathcal{A}: H^{-\sfrac12}(\partial \Omega) \to \RR$, defined via~$\mathcal{A}(\sigma\cdot \nu_\Omega) := \int_{\partial \Omega} \sigma \cdot \nu_\Omega\,d\sh^2. $ Obviously, we have 
	\begin{equation} \label{e.calA-bound}
		\bigl|\mathcal{A}(\sigma \cdot \nu_\Omega )\bigr| \leqslant \|\indc_{\partial \Omega}\|_{H^{\sfrac12}(\partial \Omega)} \|\sigma \cdot \nu_\Omega\|_{H^{-\sfrac12}(\partial \Omega)}\,. 
	\end{equation}

	\medskip \textbf{Step 2.} \textbf{The base case of the inductive construction.}
	The function~$h_0 = h + c_0 \Psi_{\ENZ}$ constructed in~\eqref{e.hconstruction-1-stmt}-\eqref{e.hconstruction-4-stmt}:
	We will define a function~$h \in H^1_c(\Omega)$ as follows. First, in the set~$D$ let~$h$ be defined to solve  
	\begin{equation} \label{e.hconstruction-1-stmt}
		\begin{dcases}
			\Delta h= 0, \quad  &\mbox{in }  D\\
			\partial_\nu h = \bF_0|_D \cdot \nu_D \quad &\mbox{ on } \partial D\,. 
		\end{dcases}
	\end{equation}
	Then, extend~$h$ in the ENZ region~$\ENZ$ by requiring that 
	\begin{equation} \label{e.hinENZ-stmt}
		\begin{dcases}
			\Delta h = 0 \quad &\mbox{ in } \ENZ\\
			h \mbox{ continuous across } &\partial D\\
			h = 0 \quad \mbox{ on }\partial \Omega\,. 
		\end{dcases}
	\end{equation}
	To complete the construction we consider the function~$\Psi_{\ENZ}$ which identically vanishes in~$D,$ and in the ENZ region~$\ENZ$ it is the unique solution to 
	\begin{equation} \label{e.hconstruction-2-stmt}
		\begin{dcases}
			-\Delta \Psi_{\ENZ} = 0 \quad &\mbox{ in } \ENZ\\
			\Psi_{\ENZ} = 0 \quad &\mbox{ on } \partial D\\
			\Psi_{\ENZ} = 1 \quad &\mbox{ on } \partial \Omega\,,
		\end{dcases}
	\end{equation}
	and observe that 
	\begin{equation} \label{e.PsiENZnonzero}
		\int_{\partial \Omega} \partial_{\nu_\Omega}\Psi_{\ENZ} \,d\sh^2 = \int_{\partial(\ENZ)} \Psi_{\ENZ} \partial_{\nu_\Omega}\Psi_{\ENZ} \,d\sh^2 = \int_{\ENZ}|\nabla \Psi_{\ENZ}|^2\,dx > 0\,.
	\end{equation}
	We choose~$c_0$ so that the vector field~$\mathbf{F}_0 + \nabla h + c_0 \nabla \Psi_{\ENZ}$ satisfies
	\begin{equation}
		\label{e.hconstruction-4-stmt}
		0 = \int_{\partial \Omega} \bF_0 \cdot \nu_\Omega +\nu_\Omega \cdot (\nabla h + c_0 \nabla \Psi_{\ENZ})\,d\sh^2 =  \int_\Omega \nu_\Omega \cdot (\nabla h + c_0 \nabla \Psi_{\ENZ})\,d\sh^2 = 0\,, 
	\end{equation}
	where we used that~$\bF_0$ is divergence free in~$\Omega,$ together with Gauss' divergence theorem. This requirement well-defines~$c_0$ in light of the computation in~\eqref{e.PsiENZnonzero}. 
	
	Thanks to~\eqref{e.hconstruction-1-stmt}-\eqref{e.hinENZ-stmt}, we obtain that 
	\begin{equation*}
		\| h\|_{H^1(D)} \leqslant C_D \|\bF_0 \cdot \nu_D\|_{H^{-\sfrac12}(\partial D)}\,,
	\end{equation*}
	and since the Neumann-to-Dirichlet map is a bounded linear operator between~$H^{-\sfrac12}(\partial D) \to H^{\sfrac12}(\partial D),$ we have 
	\begin{equation*}
		\| h\|_{H^1(\ENZ)} \leqslant C_{\ENZ} \|h\|_{H^{\sfrac12}(\partial (\ENZ))} \leqslant C_{\ENZ}C_D  \|\bF_0 \cdot \nu_D\|_{H^{-\sfrac12}(\partial D)}\,. 
	\end{equation*}
	In addition, 
	\begin{multline*}
		\biggl(\int_{\ENZ} |\nabla \Psi_{\ENZ}|^2 \biggr)   |c_0| 
		\leqslant \|\partial_{\nu_\Omega} h\|_{H^{-\sfrac12}(\partial \Omega)} \\ \leqslant C_{\ENZ} \|\nabla h\|_{H(\mathrm{div};\ENZ)} = C_{\ENZ}\|h\|_{L^2(\ENZ)} \leqslant C_{\ENZ}C_D  \|\bF_0 \cdot \nu_D\|_{H^{-\sfrac12}(\partial D)}\,,
	\end{multline*}
	since~$\nabla h$ is divergence free in~$\ENZ$ (because~$h$ is harmonic). It follows therefore that 
	\begin{equation} \label{e.basecasebound}
		\|h_0\|_{H^1(\Omega)} \leqslant C_{D,\ENZ} \|\bF_0 \cdot \nu_D\|_{H^{-\sfrac12}(\partial D)} =: C_1\,.
	\end{equation}
	This completes the base case. 
	

	In light of the operators introduced above, the objects we constructed in the base case are given by 
	\begin{equation*}
		h =
		\begin{dcases}
			T_D(\bF_0 \cdot \nu_D), \quad  \mbox{ in } D\\
			T_{\ENZ}\gamma_{0,D} T_D (\bF_0 \cdot \nu_D), \quad \mbox{ in } \ENZ\\
		\end{dcases}\,,
	\end{equation*}
	\begin{equation*}
		c_0 := - \biggl( \int_{\ENZ} |\nabla \Psi_{\ENZ}|^2\biggr)^{-1} \int_{\partial \Omega} \nu_\Omega\cdot \nabla h\,d\sh^2 = - \biggl( \int_{\ENZ} |\nabla \Psi_{\ENZ}|^2\biggr)^{-1} \mathcal{A}\bigl( \nu_\Omega \cdot \nabla  h\bigr)\,. 
	\end{equation*}
	
	\medskip
	\textbf{Step 3: Induction hypothesis.}
	We are ready to inductively construct~$\{h_j\}_{j\in \NN_0}.$ For ease of notation we allow ourselves the convention that quantities with a negative index are identically zero. 
	
	In what follows, for any function~$f$ let us denote by~$f|_D$ its restriction to~$D$ and by~$f|_{\ENZ}$ its restriction to~$\ENZ.$
	The base case of this inductive construction is~$h_0,$ and we estimated its size in~\eqref{e.basecasebound}.
	For the {induction hypothesis}, suppose we have constructed functions~$\{h_0,\cdots , h_{K-1}\}$ for some~$K \in \NN$ that satisfy 
	\begin{equation*}
		\begin{dcases}
			\mathrm{div}\, (\nabla h_k + \bF_{k}) = 0 \quad &\mbox{ in } D\\
			\nu_D \cdot (\nabla h_k + \bF_{k}) |_D= \nu_D \cdot (\nabla h_{k-1} + \bF_{k-1}\bigr)|_{\ENZ} \quad &\mbox{ on } \partial D\\
			\mathrm{div}\, (\nabla h_k + \bF_{k}) = 0, \quad &\mbox{ in } \ENZ\\
			h_k|_{\ENZ} = h_k|_D \quad &\mbox{ on }\partial D\\
			h_k = \mbox{const. } = c_k \quad &\mbox{ on }\partial \Omega\,,
		\end{dcases}
	\end{equation*}
	along with the constants~$c_k$ chosen so that 
	\begin{equation*}
		\int_{\partial \Omega} \nu_\Omega \cdot\bigl(\nabla h_k + \bF_k\bigr)\,d\sh^2 = 0, \quad \mbox{for } k = 1,\cdots, K-1\,,
	\end{equation*}
	and the functions~$\{h_k\}_{k=1}^{K-1}$ satisfy
	\begin{equation*}
		\int_D h_k \,dx = 0 \,.
	\end{equation*}
	
	\medskip
	\textbf{Step 4: Induction step.} We claim that we can solve the preceding system for~$k=K.$ 
	The choice of~$c_{K-1}$ in the induction hypothesis, and the divergence-free character of the vector field~$\bF_k$ in~$D$, render the Neumann problem in~$D$ consistent corresponding to~$k=K.$ Using the operators introduced in Step 1, we can write  
	\begin{equation} \label{e.hKdef}
		\mathring{h}_K = T_D \bigl( \nu_D \cdot \nabla h_{K-1}|_{\ENZ} + \nu_D\cdot \bF_{K-1}|_{\ENZ} - \nu_D\cdot \bF_K|_D\bigr)\,. 
	\end{equation}
	We note that this choice of a mean is simply a normalization, and it will be clear that the solution~$h_K$ is independent of this choice-- making a different choice of the average of~$h_K|D$ will simply shift the constant~$c_K$ by an amount equal to this choice. 
		Then, we let~$\mathring{h}_K|_{\ENZ}$ denote the unique solution to the problem
		\begin{equation*}
			\begin{dcases}
				\mathrm{div}\, (\nabla \mathring{h}_K + \bF_{K} )= 0 \, \quad &\mbox{ in } \ENZ\\
				\mathring{h}|_K  = \gamma_{0,D} \bigl(\mathring{h}_K |_D\bigr) = \gamma_{0,D} T_D \bigl( \nu_D \cdot \nabla h_{K-1}|_{\ENZ} + \nu_D\cdot \bF_{K-1}|_{\ENZ} - \nu_D\cdot \bF_K|_D\bigr) \, \quad &\mbox{ on } \partial D\\
				\mathring{h}|_K = 0, \quad &\mbox{ on } \partial \Omega\,,
			\end{dcases}
		\end{equation*}
		In terms of the operators introduced in Step 1, therefore we have 
		\begin{equation} \label{e.hKdefENZ}
			\mathring{h}_K|_{\ENZ} = T_{\ENZ} \gamma_{0,D} T_D \bigl( \nu_D \cdot \nabla h_{K-1}|_{\ENZ} + \nu_D\cdot \bF_{K-1}|_{\ENZ} - \nu_D\cdot \bF_K|_D\bigr) \,.
		\end{equation}
		Finally, we define~$h_K := \mathring{h}_K + c_K \Psi_{\ENZ}$ where~$c_K$ is the unique constant that verifies 
		\begin{equation*}
			\int_{\partial \Omega} \nu_\Omega \cdot \bigl( \nabla h_K + \bF_K\bigr)\,d\sh^2 = 0\,. 
		\end{equation*}
		Explicitly, this requires
		\begin{multline*}
			\int_{\partial \Omega} \partial_{\nu_\Omega} \mathring{h}_K + c_K \int_{\partial\Omega} \partial_{\nu_\Omega}\Psi_{\ENZ} + \int_{\partial\Omega}\bF_{K}\cdot \nu_\Omega=0 \\ \implies c_K = -\biggl( \int_{\ENZ}|\nabla \Psi_{\ENZ}|^2\biggr)^{-1}\biggl(\int_{\partial\Omega} (\partial_{\nu_\Omega} \mathring{h}_K + \nu_\Omega \cdot \bF_{K} )\,d\sh^2\biggr) \\
			= -\biggl( \int_{\ENZ}|\nabla \Psi_{\ENZ}|^2\biggr)^{-1}\biggl(\int_{\partial D} (\partial_{\nu_D} \mathring{h}_K + \nu_D \cdot \bF_{K} )\,d\sh^2\biggr)\,,
		\end{multline*}
		or equivalently,
		\begin{equation*}
			c_K := -\biggl( \int_{\ENZ}|\nabla \Psi_{\ENZ}|^2\biggr)^{-1}\mathcal{A} \Bigl( \partial_{\nu_D} \mathring{h}_K + \nu_D \cdot \bF_{K} \Bigr) \,. 
		\end{equation*}
		This completes the inductive construction of~$\{h_k\}_{k\in \NN}.$ 
		
		\medskip 
		\textbf{Step 5. Size of the constructed Taylor coefficients.} Towards estimating the size of the objects we constructed, we derive, by induction, a recursive inequality. As a preliminary step, let us note that since the vector field~$\nabla \mathring{h}_K + \bF_K$ is divergence-free, it has a normal trace that is in~$H^{-\sfrac12},$ and so we can estimate 
		\begin{multline*}
			|c_K| \leqslant \biggl( \int_{\ENZ}|\nabla \Psi_{\ENZ}|^2\biggr)^{-1} \biggl|\mathcal{A} \Bigl( \partial_{\nu_D} \mathring{h}_K + \nu_D \cdot \bF_{K} \Bigr) \biggr| \\
			\leqslant C_2 \|\indc_{\partial D}\|_{H^{\sfrac12}(\partial D)} \bigl\|\nu_D \cdot \bigl(\nabla \mathring{h}_K + \bF_K\bigr) \bigr) \bigr\|_{H^{-\sfrac12}(\partial D)} \\
			\leqslant C_2 \|\indc_{\partial D}\|_{H^{\sfrac12}(\partial D)} \Bigl( \|\nabla \mathring{h}_K\|_{H^1(D)} + \|\bF_K\|_{H^1(D)} \Bigr)\,. 
		\end{multline*}
		We estimate, recalling that~$\bF_K$ and~$\nabla h_{K-1}$ are divergence free in~$\Omega$ and~$\ENZ$ respectively, from~\eqref{e.hKdef}-\eqref{e.hKdefENZ} and using the properties of the operators introduced in Step 1, we obtain 
		\begin{equation*}
			\begin{aligned}
				& \|h_K\|_{H^1(\Omega)} = \|h_K\|_{H^1(\ENZ)} + \|h_K\|_{H^1(D)}\\
				&\quad  \leqslant \|\mathring{h}_K\|_{H^1(\ENZ)} + |c_K|\|\Psi_{\ENZ}\|_{H^1(\ENZ)} + \|\mathring{h}_K\|_{H^1(D)}\\
				&\quad \leqslant \|T_{\ENZ}\| \Bigl\|\gamma_{0,D} T_D \bigl( \nu_D \cdot \nabla h_{K-1}|_{\ENZ} + \nu_D \cdot (\bF_{K-1} - \bF_K)\bigr) \Bigr\|_{H^{\sfrac12}(\partial D)} \\
				&\quad \quad    + \Bigl( \|T_D\| + C_2\|\nabla \Psi_{\ENZ}\|_{H^1(\ENZ)} \|\indc_{\partial D}\|_{H^{\sfrac12}(\partial D)}\Bigr) \Bigl\| \nu_D \cdot \Bigl( \nabla h_{K-1}|_{\ENZ} + \bF_{K-1} - \bF_K \Bigr) \Bigr\|_{H^{-\sfrac12}(\partial D)}  \\
				&\quad \leqslant C_3 \Bigl\|\nu_D \cdot  \nabla h_{K-1}|_{\ENZ}\Bigr\|_{H^{-\sfrac12}(\partial D)} + C_3\Bigl\|\nu_D\cdot \Bigl(\bF_{K-1} - \bF_K\Bigr) \Bigr\|_{H^{-\sfrac12}(\partial D)}\\
				&\quad \leqslant C_3\Bigl( \Bigl\| h_{K-1}\Bigr\|_{H^1(\ENZ)} + \|\bF_{K-1}\|_{L^2(D)} + \|\bF_K\|_{L^2(D)}\Bigr)\,,
			\end{aligned}
		\end{equation*}
		where~$C_3 = \|T_{\ENZ}\|\|\gamma_{0,D}\|\|T_D\| + \|T_D\| + C_2\|\nabla \Psi_{\ENZ}\|_{H^1(\ENZ)} \|\indc_{\partial D}\|_{H^{\sfrac12}(\partial D)}. $
		
		\smallskip
		By induction, it follows that 
		\begin{equation}  \label{e.finalboundfrominduction}
			\|h_K\|_{H^1(\Omega)} \leqslant C_3^K \|h_0\|_{H^1(\Omega)} + \sum_{j=1}^{K} C_3^{1 + K-j} \Bigl( \|\bF_{j-1}\|_{L^2(D)} + \|\bF_j\|_{L^2(D)}\Bigr)\,. 
		\end{equation}

		\medskip
		\textbf{Step 6: Conclusion.} By analytic dependence of~$\bF_\de$ on~$\de,$ we know that for~$\de \in \CC, |\de|< \de_0,$ the power series~$\sum_{k=1}^\infty |\de|^k \|\bF_k\|_{L^2(\Omega)} < \infty.$ In preparation for the use of Lemma~\ref{series} below, if we set~$b_k := \|\bF_{k-1}\|_{L^2(D)} + \|\bF_{k}\|_{L^2(D)},$ we see that the series~$\sum_{k=1}^\infty |\de|^k b_k \leqslant (1 + |\de|) \sum_{k=0}^\infty \de^k \|\bF_k\|_{L^2(\Omega)} < \infty$.  
		
		\medskip
		By Lemma~\ref{series} and~\eqref{e.finalboundfrominduction},  it follows that~$\sum_{k=1}^\infty \de^k h_k$ has absolutely convergent power series for~$|\de|\leqslant \frac{\de_0}{\max(1,C_3)}. $ The proof of Lemma~\ref{l.analyticity-poisson-highcontrast} is complete. 
	\end{proof}
	In the above inductive construction we used the following lemma. 
	\begin{lemma}
		\label{series}
		Suppose~$\{b_k\}_{k\in \NN}$ is a sequence of positive numbers such that the power series
		\begin{equation*}
			f(\delta) := \sum_{k=1}^\infty b_k\delta^k\,,
		\end{equation*}
		is absolutely convergent for~$|\de|\leqslant \de_0<1.$ Then, if~$\{a_k\}_{k\in \NN}$ is a sequence of positive numbers satisfying 
		\begin{equation*}
			\begin{dcases}
				a_1 \leqslant b_1\\
				a_k \leqslant C^k b_1 + C^{k-1}b_2 + \cdots + C b_k, \quad k \in \NN\,,
			\end{dcases}
		\end{equation*}
		for some positive constant~$C.$ 
		Then the power series given by 
		\begin{equation*}
			g(\de) := \sum_{k=1}^\infty a_k \delta^k\,,
		\end{equation*}
		also converges absolutely for every~$|\delta|\leqslant \frac{\delta_0}{\max(1,C)}.$
	\end{lemma}
	\begin{proof}
		For any~$N \in \NN,$ we have, for any~$\de$ with $C|\de|< 1$ and~$|\de|\leqslant \de_0,$
		\begin{multline*}
			\Bigl| \sum_{k=1}^N a_k \delta^k \Bigr| \leqslant  \sum_{k=1}^N |a_k| |\delta|^k \leqslant \sum_{k=1}^N \Bigl( \sum_{j=1}^kC^{k-j+1} b_j\Bigr)|\delta|^k = \sum_{j=1}^N b_j\sum_{k=j}^N  C^{k-j+1}|\delta|^k \\= \frac{1}{1 - |C\de|} \sum_{j=1}^N b_j C^{-j+1} \Bigl( (C|\de|)^{j+1} - (C|\de|)^{N+1} \Bigr) 
			\leqslant \frac{C}{1-|C\de|} \sum_{j=1}^N b_j |\de|^j\,. 
		\end{multline*}
		In particular, for any~$\de$ with~$C|\de|<1$ and~$|\de|\leqslant  \delta_0,$ and any~$N \in \NN,$ the previous estimate yields
		\begin{equation*}
			\sum_{k=1}^N |a_k| |\delta|^k \leqslant \sum_{k=1}^N |a_k||\delta_0|^k \leqslant \frac{C}{1 - C\delta_0} \sum_{j=1}^\infty b_j ( \delta_0)^j\,. 
		\end{equation*}
		It follows that the series~$\sum_{k=1}^\infty a_k \delta^k$ converges absolutely for any~$\delta$ with~$|\delta|\leqslant \frac{\delta_0}{\max(1,C)}.$ 
	\end{proof}
	
	The following corollary of Lemma~\ref{l.analyticity-poisson-highcontrast} will also be used. 
	\begin{corollary} \label{c.hdel}
		Let~$\bF \in H_0(\mathrm{curl}, \Omega)$ be a divergence free vector field. Then, for~$\de \in \CC, |\de| \ll 1,$ there exist a family of functions~$h_\de \in H^1_c(\Omega)$ solving 
		\begin{equation*}
			\begin{dcases}
				\mathrm{div}\, \e_\de \bigl( \bF + \nabla h_\de\bigr) = 0, \quad \mbox{ in } \Omega\\
				h_\de = \mbox{const} \quad \mbox{ on } \partial \Omega\\
				\int_{\partial \Omega} (\bF + \nabla h_\de)\cdot \nu_\Omega\,d\sh^2 =0\,.
			\end{dcases}
		\end{equation*}
		The function~$\de \mapsto h_\de$ depends analytically on~$\de.$ 
	\end{corollary}
	\begin{proof}
		This lemma is a special case of Lemma~\ref{l.analyticity-poisson-highcontrast} when~$\bF_\de \equiv \bF$ for all~$\de.$  
	\end{proof}
	\begin{definition}
		\label{d.ellE}
		For any~$\bF \in H^1(\Omega)$ which is such that~$\int_\Omega \bF \cdot \nabla h \,dx = 0$ for all~$h \in H^1_c(\Omega),$ let~$\ell_\de(\bF) := h_\de(\bF)$ be such that 
		\begin{equation*}
			\begin{dcases}
				\mathrm{div}\, \e_\de \bigl( \bF + \nabla h_\de(\bF)\bigr) = 0, \quad \mbox{ in } \Omega\\
				h_\de(\bF) = \mbox{const} \quad \mbox{ on } \partial \Omega\,\\
				\int_{\partial \Omega} (\bF + \nabla h_\de)\cdot \nu_\Omega\,d\sh^2 =0\,.
			\end{dcases}
		\end{equation*}
		The existence and analytic dependence of the map~$\de \mapsto \ell_\de$ is a consequence of Corollary~\ref{c.hdel}. 
	\end{definition}
	\subsection{Analysis of the limiting generalized spectral problem~\eqref{e.vectorde=0}}
	In this section we solve the limiting problem~\eqref{e.vectorde=0}. Our treatment is parallel to the scalar problem, and so as in that case, we begin with functional analytic preliminaries in section~\ref{ss:FAvec}. Then, in section~\ref{ss:limviaK0} we analyze the limiting problem~\eqref{e.vectorde=0} through an intermediate operator~$\mathbb{K}_0$ that we introduce. Finally, this section concludes with perturbation theory for~\eqref{e.maxwell-delta} by another appeal to Rellich's theorem. The perturbation theory is carried out in Section~\ref{s.perturb-vector}.

\subsubsection{Functional analytic preliminaries for the Maxwell problem} \label{ss:FAvec}
This subsection is quite parallel to its scalar counterpart in~\ref{ss:FAprelim}. For any vector field~$f \in L^2(\Omega)^3$ with 
\begin{equation*}
	\int_\Omega f \cdot \nabla h \,dx = 0, \quad \mbox{ for all } h \in H^1_c(\Omega)\,,
\end{equation*}
there exists a solution~$u \in H^1(\Omega)$ to the problem 
\begin{equation*}
	\begin{dcases}
		\nabla \times \nabla \times u = f \quad \mbox{ in } \Omega\\
		u \times \nu_\Omega = 0 \quad \mbox{ on } \partial \Omega\,,
	\end{dcases}
\end{equation*}
which is unique when we impose the condition
\begin{equation*}
	\int_\Omega u \cdot \nabla h\,dx = 0 \quad \mbox{ for all } h \in H^1_c(\Omega)\,. 
\end{equation*}
It follows that the map~$f \mapsto u := (\nabla \times \nabla \times)^{-1} f$ is well-defined as a map from 
\begin{equation*}
	W := \bigl\{ \nabla h : h \in H^1_c(\Omega)\bigr\}^\perp \,,
\end{equation*}
to itself \footnote{We note that~$W$ is a closed subspace of~$L^2(\Omega;\R^3),$ and that the orthogonal complement~$\perp$ is with respect to the~$L^2$ inner product. } 
It is known that the map which takes~$f \mapsto (\nabla \times \nabla \times)^{-1} f$ is a compact operator on~$L^2(\Omega)^3$, see~\cite{Monk}; in addition it is easy to see by integration by parts that it is also self-adjoint and nonnegative. Therefore by the spectral theory for compact self-adjoint operators, the operator~$(\nabla \times \nabla \times)^{-1}$ as defined above has a sequence of eigenvalues~$\{\alpha_j\}_{j\in \NN} \in (0,\infty)$ along with associated eigenfunctions~$\{\bE_j:j \in \NN \}$ that form an orthonormal basis of~$W \subset L^2(\Omega)^3.$ The sequence~$\alpha_j \to 0^+$ as~$j\to \infty,$ and is arranged in nonincreasing order, repeating according to finite multiplicity of multiple eigenvalues. Moreoever, the vectors~$\{\bE_j\}_{j\in \NN}$ are also orthogonal with respect to the inner product 
\begin{equation} \label{e.Hcurlinnerp}
	\langle E, F \rangle := \int_\Omega \bigl( \nabla \times E \bigr)\cdot \bigl( \nabla \times F \bigr)\,dx + \int_\Omega (\mathrm{div}\, E)(\mathrm{div}\, F)\,dx\,. 
\end{equation}
It can be checked by integration by parts that for vector fields~$E,F$ with vanishing tangential component along~$\partial \Omega,$ this inner product is equivalent to the standard~$H^1(\Omega)^3$ inner product, and 
\begin{equation*}
	\langle \bE_j, \bE_j\rangle = \int_\Omega |\nabla \times \bE_j|^2\,dx = \alpha_j \int_\Omega |\bE_j|^2\,dx = \alpha_j\,.  
\end{equation*}
Therefore, the set 
\begin{equation*}
	\Bigl\{ \frac{1}{\sqrt{\alpha_j}} \bE_j \Bigr\}_{j\in \NN}
\end{equation*}
forms an orthonormal set with respect to the inner product in~\eqref{e.Hcurlinnerp}. Therefore, as in Section~\ref{ss:FAprelim}, we may define fractional powers of the operator~$(\nabla \times \nabla \times)$ spectrally: namely, we define 
\begin{multline} \label{e.Vdef}
	\bigl( \nabla \times \nabla \times )^{-\sfrac12} : W \to V \,, \\
	V:= \Bigl\{ E \in H^1(\Omega) : E \times \nu_\Omega = 0 \quad \mbox{ at } \partial \Omega, \int_\Omega E \cdot \nabla h \,dx = 0 \mbox{ for all } h \in H^1_c(\Omega)\Bigr\}\,,
\end{multline}
by setting, for any~$E \in W$ with~$E = \sum_{j=1}^\infty c_j \bE_j,$
\begin{equation*}
	\bigl( \nabla \times \nabla \times\bigr)^{-\sfrac12} E := \sum_{j=1}^\infty \frac{c_j}{\sqrt{\alpha_j}} \bE_j\,. 
\end{equation*}
This map is surjective, and we can also define 
\begin{equation*}
	\bigl( \nabla \times \nabla \times)^{\sfrac12} : V \to W\,,
\end{equation*}
by setting: 
\begin{equation*}
	\bigl( \nabla \times \nabla \times)^{\sfrac12} E := \sum_{j=1}^\infty \sqrt{\alpha_j} c_j \bE_j\,.
\end{equation*}
Parallel to Lemma~\ref{l.observations},  these maps are isometries between~$V$ with the inner product~$\langle \cdot, \cdot \rangle$ from~\eqref{e.Hcurlinnerp}, and~$W$ with the~$L^2(\Omega)^3$ inner product. Finally, parallel to Lemma~\ref{l.equivalence-2}, we have:
\begin{lemma}
	\label{l.maxwellequivalence}
	If~$E \in W$ and~$F \in V$ are such that 
	\begin{equation*}
		\bigl( \nabla \times \nabla \times)^{-\sfrac12} E = (\nabla \times \nabla \times)^{\sfrac12} F\,,
	\end{equation*}
	then in fact, 
	\begin{equation}
		\begin{dcases}
			\nabla \times \nabla \times F = E \quad \mbox{ in } \Omega\\
			F \times \nu_\Omega = 0 \quad \mbox{ on } \partial \Omega\,,
		\end{dcases}
	\end{equation}
	in the sense that 
	\begin{equation*}
		\int_\Omega \nabla \times F \cdot \nabla \times G = \int_\Omega E \cdot G\,, \quad \mbox{ for all } G \in H^1(\Omega) \mbox{ such that } G \times \nu_\Omega = 0 \quad \mbox{on } \Omega\,. 
	\end{equation*}
\end{lemma}
\begin{proof}
	The proof of this lemma is the same as that of Lemma~\ref{l.equivalence-2}, up to notational changes, and is therefore omitted.
\end{proof}
\smallskip 
\subsubsection{Analysis of~\eqref{e.vectorde=0} via an operator~$\mathbb{K}_0$} \label{ss:limviaK0}
We define the linear operator~$\mathbb{K}_0:W \to W$ via
\begin{equation} \label{e.K0defvec}
	\mathbb{K}_0 := (\nabla \times \nabla \times)^{-\sfrac12}(\P_\Omega \indc_D) (\nabla \times \nabla \times)^{-\sfrac12}\,. 
\end{equation}
This operator extends to a compact operator from~$L^2(\Omega)^3$ to itself, if we set
\begin{equation*}
	\mathbb{K}_0 = \P_\Omega^T (\nabla \times \nabla \times)^{-\sfrac12} \P_\Omega \indc_D \P_\Omega^T (\nabla \times \nabla \times)^{-\sfrac12} \P_\Omega\,,
\end{equation*}
where, we recall, from~\eqref{e.Pomegadef}, that~$\P_\Omega: L^2(\Omega)^3 \to W$ is orthogonal projection. 

It is clear from the definition that~$\mathbb{K}_0$ is a self-adjoint operator on~$L^2(\Omega)^3;$ Furthermore, since~$\P_\Omega \indc_D$ is a bounded linear operator and~$(\nabla \times \nabla \times)^{-\sfrac12}$ is a compact operator, it follows that the operator~$\mathbb{K}_0$ is compact.  we check its positivity in the next lemma. 
\begin{lemma}
	\label{l.bbkernel}
	The operator~$\mathbb{K}_0$ defined above is a nonnegative operator in the sense that for every~$v \in L^2(\Omega)^3$ we have 
	\begin{equation*}
		\int_\Omega (\mathbb{K}_0 v ) v \,dx \geqslant 0\,.
	\end{equation*}
	In addition, the kernel of~$\mathbb{K}_0$ is given by
	\begin{equation*}
		\ker (\mathbb{K}_0) = \biggl\{ \bigl( \nabla \times \nabla \times)^{\sfrac12} \xi + \nabla h: h \in H^1_c(\Omega) \quad 
		\mbox{ and } \xi \in V : \xi \equiv 0 \mbox{ on } D\biggr\}\,.
	\end{equation*}
\end{lemma}
\begin{proof}
	To check positivity of~$\mathbb{K}_0$ we may restrict to~$W$ since the definition of~$\mathbb{K}_0$ begins by projection to~$W.$ Since~$\{\bE_j\}_{j\in \NN}$  form an orthonormal basis of~$W$, writing~$v$ as 
	\begin{equation*}
		v = \sum_{j=1}^\infty \biggl( \int_\Omega v \cdot \bE_j \,dx\biggr)\bE_j\,.
	\end{equation*}
	Therefore, 
	\begin{equation*}
		\indc_D \P_\Omega^T \bigl( \nabla \times \nabla \times \bigr)^{-\frac12}v = \sum_{j=1}^\infty \frac1{\sqrt{\alpha_j}} \biggl(\int_\Omega v \cdot \bE_j \biggr) \bE_j \indc_D\,. 
	\end{equation*}
	As discussed at the start of Section~\ref{s.helmholtz}, the projection of this vector field to~$W$ has the form 
	\begin{equation*}
		\Xi :=   \sum_{j=1}^\infty \frac1{\sqrt{\alpha_j}} \biggl( \int_\Omega v \cdot \bE_j\biggr) \bE_j \indc_D - \nabla q\,,
	\end{equation*}
	with~$q$ chosen so that 
	\begin{equation*}
		\Xi \perp \bigl\{ \nabla h : h \in H^1_c(\Omega)\bigr\}\,.
	\end{equation*}
	Now, we compute, using that each~$\bE_j \perp \nabla H^1_c(\Omega),$ that 
	\begin{align*}
		\int_\Omega \mathbb{K}_0 v \cdot v \,dx &= \int_\Omega \biggl( \sum_{j=1}^\infty \frac1{\sqrt{\alpha_j}} \biggl( \int_\Omega v \cdot \bE_j \,dx\biggr)\bE_j \indc_D- \nabla q \biggr)\cdot \biggl( \sum_{j=1}^\infty \frac1{\sqrt{\alpha_j}} \biggl(\int_\Omega v \cdot \bE_j \,dx\biggr)\bE_j\biggr)\\
		& =  \int_\Omega \biggl( \sum_{j=1}^\infty \frac1{\sqrt{\alpha_j}} \biggl( \int_\Omega v \cdot \bE_j \biggr)\bE_j \indc_D\biggr) \cdot \biggl( \sum_{j=1}^\infty \frac1{\sqrt{\alpha_j}} \biggl(\int_\Omega v \cdot \bE_j \,dx\biggr)\bE_j\biggr)\\
		&= \int_D \biggl|\sum_{j=1}^\infty \frac1{\sqrt{\alpha_j}} \biggl(\int_\Omega v \cdot \bE_j\,dx\biggr)\bE_j \biggr|^2\geqslant 0\,.
	\end{align*}
	To characterize the kernel of~$\mathbb{K}_0:$ from the preceding computation, 
	\begin{equation*}
		\mathbb{K}_0 v = 0 \implies \ \ \bigl( \nabla \times \nabla \times\bigr)^{-\sfrac12} v = 0 \quad \mbox{ in } D\,. 
	\end{equation*}
	Setting~$\xi = \bigl( \nabla \times \nabla \times\bigr)^{-\sfrac12} v,$ we infer that~$v = \bigl( \nabla \times \nabla \times\bigr)^{\sfrac12} \xi$ for some~$\xi \in V$ with~$\xi \equiv 0$ on~$D.$ 
	
	Conversely, if~$v \in W$ satisfies ~$v = \bigl( \nabla \times \nabla \times\bigr)^{\sfrac12} \xi$ for some~$\xi \in V$ with~$\xi \equiv 0$ on~$D,$ then 
	\begin{equation*}
		\indc_D \P_\Omega^T \bigl( \nabla \times \nabla \times )^{-\sfrac12} \P_\Omega v = \indc_D \xi = 0\,.
	\end{equation*}
	This implies that~$\mathbb{K}_0 v = 0.$ 
	
	Finally, since~$\mathbb{K}_0$ begins with projection to~$W,$ as an operator on~$L^2(\Omega)^3,$ we have 
	\begin{equation*}
		\ker (\mathbb{K}_0) = \biggl\{ \bigl( \nabla \times \nabla \times)^{\sfrac12} \xi + \nabla h: h \in H^1_c(\Omega) \quad 
		\mbox{ and } \xi \in V : \xi \equiv 0 \mbox{ on } D\biggr\}\,.
	\end{equation*}
	The proof is completed. 
\end{proof}

Thanks to the spectral theorem of compact self-adjoint operators, the spectrum of the operator~$\mathbb{K}_0$ consists of a sequence of nonzero eigenvalues that are positive, say~$\{\rho_{0,k}\}_{k\in \NN} $, each with finite multiplicity, along with associated eigenfunctions~$\{\chi_{0,k}\}_{k\in \NN}$ that form an orthonormal Hilbert basis of~$\ker(\mathbb{K}_0)^\perp.$

As in our scalar problem, we can use these eigenvalues to construct eigenvalues and eigenfunctions for~\eqref{e.vectorde=0}. To see how, let us note that if~$\bF_0$ is an eigenfunction of the compact operator~$\mathbb{K}_0$ associated to eigenvalue~$\la^{-1},$ then define~$\tilde{\bE}_0 :=   (\nabla \times \nabla \times)^{-\sfrac12} \bF_0,$ and note that this vector field is divergence free. 

We will define a function~$h \in H^1_c(\Omega)$ as follows. First, in the set~$D$ let~$h$ be defined to solve  
\begin{equation} \label{e.hconstruction-1}
	\begin{dcases}
		\Delta h= 0, \quad  &\mbox{in }  D\\
		\partial_\nu h = \tilde{\bE}_0|_D \cdot \nu_D \quad &\mbox{ on } \partial D\,. 
	\end{dcases}
\end{equation}
Then, define~$h$ in the ENZ region~$\ENZ$ by requiring that 
\begin{equation} \label{e.hinENZ}
	\begin{dcases}
		\Delta h = 0 \quad &\mbox{ in } \ENZ\\
		h \mbox{ continuous across } &\partial D\\
		h = 0 \quad \mbox{ on }\partial \Omega\,. 
	\end{dcases}
\end{equation}
To complete the construction we consider the function~$\Psi_{\ENZ}$ which identically vanishes in~$D,$ and in the ENZ region~$\ENZ$ it is the unique solution to 
\begin{equation} \label{e.hconstruction-2}
	\begin{dcases}
		-\Delta \Psi_{\ENZ} = 0 \quad &\mbox{ in } \ENZ\\
		\Psi_{\ENZ} = 0 \quad &\mbox{ on } \partial D\\
		\Psi_{\ENZ} = 1 \quad &\mbox{ on } \partial \Omega\,,
	\end{dcases}
\end{equation}
and observe that 
\begin{equation*}
	\int_{\partial \Omega} \partial_{\nu_\Omega}\Psi_{\ENZ} \,d\sh^2 = \int_{\partial(\ENZ)} \Psi_{\ENZ} \partial_{\nu_\Omega}\Psi_{\ENZ} \,d\sh^2 = \int_{\ENZ}|\nabla \Psi_{\ENZ}|^2\,dx > 0\,.
\end{equation*}
We then set~$\bE = \tilde{\bE}_0 - \nabla h - c_0 \nabla \Psi_{\ENZ},$ where the constant~$c_0$ is chosen to arrange that
\begin{equation}
	\label{e.hconstruction-4}
	\int_{\partial \Omega} \bE \cdot \nu\,d\sh^2 = 0 = \int_{\partial \Omega} \tilde{\bE}_0 \cdot \nu_\Omega - \nu_\Omega \cdot (\nabla h + c_0 \nabla \Psi_{\ENZ})\,d\sh^2 = - \int_\Omega \nu_\Omega \cdot (\nabla h + c_0 \nabla \Psi_{\ENZ})\,d\sh^2 = 0\,, 
\end{equation}
where in the second equality we used that~$\tilde{\bE}_0$ is divergence free in~$\Omega,$ together with Gauss' divergence theorem. This requirement well-defines~$c_0$ in light of the computation in the previous line. 
With this definition, we have~$\mathrm{div}(\bE \indc_D) = 0$ in~$\Omega,$ or equivalently,~$\P_\Omega \indc_D(\tilde{\bE}_0 - \nabla h) = \indc_D \bE$, and can compute 
\begin{equation*}
	\nabla \times \nabla \times \bE = \nabla \times \nabla \times \tilde{\bE}_0  = (\nabla \times \nabla \times)^{\sfrac12} \bF_0 = \la \P_\Omega \indc_D (\nabla \times \nabla \times)^{-\sfrac12} \bF_0 = \la \P_\Omega \indc_D \tilde{\bE}_0 = \la \indc_D \bE\,,
\end{equation*}
in~$\Omega$, along with~$\bE \times \nu_\Omega = 0$ on~$\partial \Omega.$ This says that~$(\bE, \la)$ is solves the generalized eigenvalue problem~\eqref{e.vectorde=0}. 

\smallskip 
Conversely, the foregoing procedure of constructing an eigenpair~$(\bE,\la)$ solving the generalized eigenproblem~\eqref{e.vectorde=0} from an eigenpair of the compact operator~$\mathbb{K}_0$ can be reversed-- and starting from an eigenpair~$(\bE,\la)$ of~\eqref{e.vectorde=0} one can construct an eigenpair for~$\mathbb{K}_0.$ This yields a bijective correspondence between the two sets of eigenpairs and this procedure captures all eigenvalues of~\eqref{e.vectorde=0}.  

We have thus shown: 
\begin{theorem}
	\label{t.maxwelllim-spectrum}
	There exist a sequence of positive numbers~$\{\la_{0,k}\},$ arranged in nondecreasing order, repeated according to multiplicity (each finite) with~$\la_{0,k} \to \infty$ as~$k \to \infty,$ and associated eigenfunctions~$\bE_{0,k} \in H_0(\mathrm{curl};\Omega)$ with~$\mathrm{div} \, (\bE_{0,k}\indc_D) = 0$ satisfying~\eqref{e.vectorde=0}. 
\end{theorem}
\begin{proof}
	We must set~$\la_{0,k} = \rho_{0,k}^{-1}$ where~$\rho_{0,k}>0$ is the~$k$th eigenvalue of the compact operator~$\mathbb{K}_0$ and let~$\bF_{0,k}$ the associated eigenfunction of the operator~$\mathbb{K}_0.$ Then one carries out the discussion prior to this present theorem statement to construct~$\bE_{0,k}.$ 
\end{proof}
\subsection{Perturbation theory for~\eqref{e.maxwell-delta}} \label{s.perturb-vector}
We turn to analyzing~\eqref{e.maxwell-delta} for~$|\de| \ll 1.$ Once again, we study this generalized eigenproblem in two steps: we first obtain a description of the spectrum of a nonlocal operator~$\mathbb{K}_\de$ given by: for any~$E \in W,$ we define
\begin{equation} \label{e.vecK0def}
	\mathbb{K}_\de(E) := (\nabla \times \nabla \times)^{-\sfrac12}  \e_\de \Bigl( (\nabla \times \nabla \times)^{-\sfrac12} E + \nabla h_\de (E) \Bigr)\,,
\end{equation}
where~$h_\de(E) := \ell_\de\bigl( (\nabla \times \nabla \times)^{-\sfrac12} E\bigr) \in H^1_c(\Omega)$ is the linear operator defined in Definition~\ref{d.ellE}. It is the unique function~$h_\de \in H^1_c(\Omega)$ which satisfies 
\begin{equation*}
	\int_\Omega \e_\de \Bigl( (\nabla \times \nabla \times)^{-\sfrac12} E + \nabla h_\de (E) \Bigr) \cdot \nabla k \,dx = 0 \quad \mbox{ for all } k \in H^1_c(\Omega)\,. 
\end{equation*}
By corollary~\ref{c.hdel}, the map~$\de \mapsto h_\de(E)$ is analytic. 

\begin{lemma}[Self-adjointness of~$\mathbb{K}_\de$]
	\label{l.Kdelselfadj}
	The operator~$\mathbb{K}_\de$ is a self adjoint operator on~$W.$ In addition it is a compact operator on~$W.$ 
\end{lemma}
\begin{proof}
	If~$E,F \in W,$ then 
	\begin{equation*}
		\begin{aligned}
			& \int_\Omega \mathbb{K}_\de (E) \cdot F\,dx \\
			&\quad = \int_\Omega (\nabla \times \nabla \times)^{-\sfrac12}  \e_\de \Bigl( (\nabla \times \nabla \times)^{-\sfrac12} E + \nabla h_\de (E) \Bigr) \cdot F\,dx \\
			&\quad = \int_\Omega  \e_\de \Bigl( (\nabla \times \nabla \times)^{-\sfrac12} E + \nabla h_\de (E) \Bigr) \cdot (\nabla \times \nabla \times)^{-\sfrac12} F \\
			&\quad = \int_\Omega  \e_\de \Bigl( (\nabla \times \nabla \times)^{-\sfrac12} E + \nabla h_\de (E) \Bigr) \cdot \Bigl(  (\nabla \times \nabla \times)^{-\sfrac12} F + \nabla h_\de (F) \Bigr)\,dx = \int_\Omega E \cdot \mathbb{K}_\de(F)\,dx \,,
		\end{aligned}
	\end{equation*}
	where in the last line we used that~$\nabla h_\de (E)$ is chosen exactly to make~$\e_\de \Bigl( (\nabla \times \nabla \times)^{-\sfrac12} E + \nabla h_\de (E) \Bigr)$ orthogonal in~$L^2$ to all gradients~$\nabla k, k \in H^1_c(\Omega),$ and, in particular,~$h_\de(F) \in H^1_c(\Omega). $ 
	
	\medskip 
	The statement about the compactness of ~$\mathbb{K}_\de$ on~$W$ is well-known, see for e.g.~\cite{Monk}. The proof of the lemma concludes.
\end{proof}
We are now ready to prove the main theorem of this section, which concerns the perturbation theory for the spectrum of the operator~$\mathbb{K}_\de,$ yielding eigenvalues and eigenfunctions of the Maxwell system~\eqref{e.strongform}.

Naturally, our strategy is quite parallel to the scalar case, and we first analyze the analytic dependence of the spectrum of the operator~$\mathbb{K}_\de.$ For the dependence of the spectrum of~$K_\de$ on~$\de,$ we note:
\begin{itemize}
	\item by Lemma~\ref{l.Kdelselfadj} this is a family of self-adjoint operators  depending analytically on~$\de$ (by Lemma~\ref{l.analdepofldel})\,,
	\item each~$\mathbb{K}_\de$ is also compact on~$W.$ To see this, by~\cite{Monk}, the operator~$(\nabla \times \nabla \times)^{-1}$ is compact. Then, for~$(\nabla \times \nabla \times)^{-\sfrac12}$, and hence for~$\mathbb{K}_\de$ one argues as in the proof of Lemma~\ref{l.K0good}. 
\end{itemize}

Rellich's theorem, recalled in~Theorem \ref{t.rellich}, immediately applies. To state the conclusions that we can draw from Rellich's theorem, let us introduce some notation. First, let~$\la_{0,k}^{-1}$ denote the~$k$th eigenvalue of~$\mathbb{K}_0$ introduced in~\eqref{e.K0defvec} and let~$\la_{\de,k}^{-1}$ denote the~$k$th eigenvalue of~$\mathbb{K}_\de$ introduced in~\eqref{e.vecK0def}. Let the associated eigenfunctions of~$\mathbb{K}_0$ be~$\bE_{0,k}$ (for~$\mathbb{K}_0$) and~$\bE_{\de,k}$ (for~$\mathbb{K}_\de)$ respectively. We know that~$\la_{0,k}^{-1}>0$ and we know from the spectral theorem for compact self-adjoint operators that each of these eigenvalues has finite multiplicity. 

For the limiting operator~$\mathbb{K}_0$, we already know from Theorem~\ref{t.maxwelllim-spectrum} that we can build eigenfunctions for~\eqref{e.vectorde=0} from an eigenfunction of~$\mathbb{K}_0$, and that all eigenfunctions of~\eqref{e.vectorde=0} are obtained this way. The proof of the next theorem will show how to obtain eigenfunctions of~\eqref{e.maxwell-delta} from eigenfunctions of~$\mathbb{K}_\de$-- the heart of the matter is the~$\de-$dependent Helmholtz projection from Lemma~\ref{l.analyticity-poisson-highcontrast} that is analytic in~$\de.$ 

Let~$\gamma(\la_{0,k})$ denote the spectral gap about~$\la_{0,k}:$
\begin{equation*}
	\ga(\la_{0,k}) := \min_{j \in \NN} \{|\la_{0,j} - \la_{0,k}| : \la_{0,j} \neq \la_{0,k}\}>0\,. 
\end{equation*}
The main theorem of this section is:
\begin{theorem}
	\label{t.vectormain}
	For each~$k \in \NN$, let~$\lambda_{0,k}$ have (finite) multiplicity~$h,$ and let~$0 < \gamma < \gamma(\la_{0,k})$ be fixed, so that~$\mathbb{K}_0$ has no eigenvalues in closed neighborhood~$[\gamma(\la_{0,k}) - \gamma, \gamma(\la_{0,k}) + \gamma]$ other than~$\la_{0,k}$, and suppose, without loss of generality, that~$\la_{0,k-1} < \la_{0,k} = \la_{0,k+1} = \cdots = \la_{0,k+h-1} < \la_{0,k+h}.$ Then there exists a positive number~$\delta_0 = \delta_0(k),$ and~$h$ families of complex numbers 
	\begin{equation*}
		\la_{\de,k}, \cdots, \la_{\de,k+h-1}\,,
	\end{equation*}
	that have convergent power series in~$\de$ in a neighborhood~$|\de|< \de_0(k)$; additionally, there exist an associated family of functions~$\{u_{\de,k},\cdots, u_{\de,k+h-1}\} \subset W$ also given by a convergent power series in the Hilbert space~$L^2(\Omega)^3$ that satisfy the following properties: 
	For each~$i = 0,\cdots, h-1,$ we have 
	\begin{equation} \label{e.udelk+i}
		\begin{dcases}
			\nabla \times \nabla \times u_{\de,k+i} &= \la_{\de,k+i} \e_\de u_{\de,k+i}, \quad \mbox{ in } \Omega\\
			\nu_\Omega \times u_{\de,k+i} &= 0, \quad \quad \mbox{ on } \partial \Omega\,.
		\end{dcases}
	\end{equation}
	Furthermore, for each~$i \in \{0,\cdots, h-1\}$, we have that~$\lim_{\de \to 0}\la_{\de,k+i} = \la_{0,k+i},$ and in addition, for all~$\de, |\de|<\de_0(k),$ and~$i,j = 0,\cdots,h-1,$ the following orthogonality relation holds:
	\begin{equation*}
		\int_{\Omega}\e_\de u_{\de,k+i} \cdot u_{\de,k+j} = \de_{ij} = \begin{dcases}
			1 & \mbox{ if } i = j\\
			0 & \mbox{ otherwise } 
		\end{dcases}\,. 
	\end{equation*}
	Finally, for~$|\de| \ll 1,$ the interval~$[\la_{0,k} - \frac{\gamma}2, \la_{0,k} + \frac{\gamma}2]$ includes no other eigenvalues for the generalized eigenvalue problem~\eqref{e.PDEscalar} except for the~$h$ branches~$\la_{\de,k}, \cdots, \la_{\de,k+h-1}.$ 
\end{theorem}
\begin{proof}
	By Lemmas~\ref{l.analyticity-poisson-highcontrast} and~\ref{l.Kdelselfadj}, the family of operators~$\de \mapsto \mathbb{K}_\de$ is a family of compact operators on~$\mathring{L}^2(\Omega),$ which are self-adjoint for~$|\de| \ll 1$ when~$\de \in \RR.$ Therefore Rellich theorem applies, and immediately yields~$h$ analytic branches~$\la_{\de,k}^{-1}, \cdots, \la_{\de,k+h-1}^{-1}$ of eigenvalues and associated eigenfunctions~$\{\bE_{\de,k}, \cdots , \bE_{\de,k+h-1}\}$ of the operator~$K_\de$ that all admit convergent power series representations in~$\de$ for~$|\de| < \de_0(k)$ sufficiently small. In addition, these satisfy~$\lim_{\de \to 0} \la_{\de,k+i} = \la_{0,k+i}$ for each~$i=0,\cdots,h-1.$ Since~$\la_{0,k+i} > 0$ for each~$i,$ it follows that the maps~$\de \mapsto \la_{\de,k}$ are also analytic. 
	
	\medskip 
	Finally, by Proposition~\ref{l.maxwellequivalence}, since for any~$r \in \{0,\cdots, h-1\},$
	\begin{equation*}
		\mathbb{K}_\de \bE_{\de,k+r} = \bigl( \nabla \times \nabla \times\bigr)^{-\sfrac12} \e_\de \Bigl( \bigl(\nabla \times \nabla \times \bigr)^{-\sfrac12} \bE_{\de,k+r} + \nabla h_\de (\bE_{\de,k+r})  \Bigr) = \la_{\de,k+r}^{-1}\bE_{\de,k+r}\,.
	\end{equation*}
	Define
	\begin{equation} \label{e.udeldef}
		u_{\de,k+r} := \bigl(\nabla \times \nabla \times \bigr)^{-\sfrac12} \bE_{\de,k+r} + \nabla h_\de (\bE_{\de,k+r}) \,.
	\end{equation}
	Applying~$(\nabla \times \nabla \times)$ on both sides, we find,
	\begin{multline*}
		\nabla \times \nabla \times u_{\de,k+r} 
		= (\nabla \times \nabla \times)^{\sfrac12} \bE_{\de,k+r} \\
		= \la_{\de,k+r} \e_\de \Bigl( \bigl(\nabla \times \nabla \times \bigr)^{-\sfrac12} \bE_{\de,k+r} + \nabla h_\de (\bE_{\de,k+r})\Bigr)  = \la_{\de,k+r} \e_\de u_{\de,k+r}\,.
	\end{multline*}
	To check the boundary conditions, we recall from Lemma~\ref{l.maxwellequivalence} and the prior discussion that~$(\nabla \times \nabla \times)^{-\sfrac12}$ maps into~$V$ (defined in~\eqref{e.Vdef}), and that~$h_\de (\bE_{\de,k+r}) \in H^1_c(\Omega)$ is constant along~$\partial \Omega,$ and so has vanishing tangential derivatives. 
	
	Therefore, 
	\begin{equation*}
		\nu_\Omega \times u_{\de,k+r} = \nu \times \Bigl(\bigl(\nabla \times \nabla \times \bigr)^{-\sfrac12} \bE_{\de,k+r} + \nabla h_\de (\bE_{\de,k+r})  \Bigr) = 0 \, \quad \mbox{ on } \partial \Omega\,. 
	\end{equation*}
	
	Finally, we must check the orthogonality condition. If~$\lambda_{\de,k+r} \neq \lambda_{\de,k+s}$ then this follows from multiplying the equation for~$u_{\de,k+r}$ by~$u_{\de,k+s}$ and vice-versa, integrating, and subtracting the results:
	\begin{equation*}
		0 = \int_{\Omega } \nabla u_{\de,k+r} \cdot \nabla u_{\de,k+s} \,dx = \bigl( \la_{\de,k+r} - \la_{\de,k+s}\bigr) \int_\Omega \e_\de u_{\de,k+r}u_{\de,k+s}\,dx\,. 
	\end{equation*}
	If, for some~$r \neq s,$ it happens that~$\la_{\de,k+r} = \la_{\de,k+s}$ also for (some)~$\de \neq 0$ (as it might, rarely), then there is a different argument to prove the same orthogonality: using~\eqref{e.udeldef}, the facts that~$\nabla h_\de (\bE_{\de,k+r})$ is a gradient, and that~$\nabla h_\de(\bE_{\de,k+r})$ is chosen to make~$\e_\de \Bigl((\nabla \times \nabla \times)^{-\sfrac12} \bE_{\de,k+r} + \nabla h_\de (\bE_{\de,k+r})\Bigr)$ orthogonal to gradients of~$H^1_c(\Omega)$ functions (and in particular,~$\nabla h_\de(\bE_{\de,k+s})$), and finally, the self-adjointness of~$(\nabla \times \nabla \times)^{-\sfrac12},$ we compute
	\begin{multline}
		\int_\Omega \e_\de u_{\de,k+r} \cdot  u_{\de,k+s}\,dx \\
		= \int_\Omega \e_\de \Bigl( \bigl(\nabla \times \nabla \times \bigr)^{-\sfrac12} \bE_{\de,k+r} + \nabla h_\de (\bE_{\de,k+r}) \Bigr) \cdot \Bigl( \bigl(\nabla \times \nabla \times \bigr)^{-\sfrac12} \bE_{\de,k+s} + \nabla h_\de (\bE_{\de,k+s}) \Bigr)\,dx \\
		= \int_\Omega \e_\de  \bigl(\nabla \times \nabla \times \bigr)^{-\sfrac12} \bE_{\de,k+r} \cdot  \Bigl( \bigl(\nabla \times \nabla \times \bigr)^{-\sfrac12} \bE_{\de,k+s} + \nabla h_\de (\bE_{\de,k+s})\Bigr)\,dx \\
		= \int_\Omega   \bigl(\nabla \times \nabla \times \bigr)^{-\sfrac12} \bE_{\de,k+r} \cdot \e_\de \Bigl( \bigl(\nabla \times \nabla \times \bigr)^{-\sfrac12} \bE_{\de,k+s} + \nabla h_\de (\bE_{\de,k+s})\Bigr)\,dx\\
		= \int_\Omega    \bE_{\de,k+r} \cdot \bigl(\nabla \times \nabla \times \bigr)^{-\sfrac12}\e_\de \Bigl( \bigl(\nabla \times \nabla \times \bigr)^{-\sfrac12} \bE_{\de,k+s} + \nabla h_\de (\bE_{\de,k+s})\Bigr)\,dx\\
		= \la_{\de,k+s}^{-1} \int_\Omega \bE_{\de,k+r} \cdot \bE_{\de,k+s}\,dx = 0\,,
	\end{multline}
	since~$\bE_{\de,k+r}$ and~$\bE_{\de,k+s}$ are~$L^2(\Omega)^3-$orthogonal by virtue of being eigenfunctions of the self-adjoint operator~$\mathbb{K}_\de.$ 
\end{proof}

\section{Electrostatic resonances} \label{sec:electro}
In this section we describe a specific class of resonances identified by~\cite{Engheta}-- \emph{electrostatic resonances}, which have the property that at leading order, the magnetic field \emph{vanishes}, and so the electric field is the gradient of a harmonic function in the ENZ region~$\ENZ$ which is constant at the outer boundary~$\partial \Omega.$ Examining~\eqref{e.strongform}, it is clear that in order for a resonance~$\la_{0,j}$ to be electrostatic to leading order, the eigenfunction~$\bE_{0,j}$ must solve an overdetermined vectorial problem in~$D:$
\begin{equation} \label{e.overdet}
	\begin{dcases}
		\nabla \times \nabla \times \bE_{0,j} = \la_{0,j} \bE_{0,j} & \quad \mbox{ in } D\\
		\bE_{0,j} \cdot \nu_D = 0 &\quad \mbox{ on } \partial D\\
		\bH_{0,j} := \nabla \times \bE_{0,j}  = 0 & \quad \mbox{ on } \partial D\,. 
	\end{dcases}
\end{equation}
We say that an inclusion~$D$ \emph{supports an electrostatic resonance} if the overdetermined problem~\eqref{e.overdet} has a nontrivial solution (i.e., there exists a nonzero~$\la_{0,j}$ such that~\eqref{e.overdet} holds). 

\begin{remark}
	\label{r.why-overdetermined}
	To see why~\eqref{e.overdet} represents an overdetermined system of equations, we recall from the literature on Maxwell's equations (see for e.g.~\cite{costabeldauge}) that the cavity eigenvalue problem for Maxwell's equations corresponds to the ``perfectly conducting'' boundary conditions~$\bH \times \nu_D = 0$ and~$\bE \cdot \nu_D = 0.$ We observe that in addition to these, the system~\eqref{e.overdet} includes the condition~$\bH \cdot \nu_D = 0,$ rendering the problem overdetermined. 
\end{remark}
\begin{remark}\label{r.electrostaticresonance}
	A natural question is to ask which domains~$D \subset \R^3$ support an electrostatic resonance. This problem is \emph{open}. What we can show below, is that if~$D$ is a ball, then it supports infinitely many electrostatic resonances, each that has high multiplicity. 
\end{remark}

Let us explain how electrostatic resonances lead to ``geometry invariant closed resonant cavities": Let~$\Omega \subset \R^3$ denote a Lipschitz domain that contains the closure of~$D.$ Then, given an electrostatic resonance, it is straightforward to extend such a resonance to a solution to~\eqref{e.strongform} by simply solving the equation 
\begin{equation*}
	\begin{dcases}
		\nabla \times \nabla \times \bE_{0,j} = 0  &\quad \mbox{ in } \ENZ\\
		\bE_{0,j} \times \nu_D \mbox{ is continuous across } &\partial D\\
		\bE_{0,j} \times \nu_\Omega = 0 &\quad \mbox{ on } \partial \Omega\,. 
	\end{dcases}
\end{equation*}
The geometry invariance refers to the statement that~$\la_{0,j}$ is \emph{independent of the shape of~$\Omega.$ } Accordingly, a consequence of Theorem~\ref{t.electro} below is that for a core-shell arrangement with a spherical core, there are infinitely many electrostatic resonances with ENZ material in the shell. 

\begin{remark}
	When~$\de \neq 0,$ by Theorem~\ref{t.vectormain}, any of the electrostatic resonances we construct in this section bifurcate for~$0 < |\de|\ll 1,$ to generate quasiresonances of~\eqref{e.strongform} that decay in time when the imaginary part of~$\de$ is negative (i.e., in the presence of losses). Which of these shapes optimizes the life time of the resonance? We explored this question in two dimensions (i.e., tranverse magnetic resonances) in~\cite{KV2}. We leave an exploration of the corresponding three dimensional question to future work. 
\end{remark}
In order to facilitate a self-contained presentation and in order to fix notation, we begin our analytical construction of electrostatic resonances with a digression on spherical harmonics. 

\smallskip
\emph{Digression on spherical harmonics:}
For each~$n \in \NN$ the Bessel function~$j_n$ solves the ordinary differential equation 
\begin{equation*}
	j_n''(r) + \frac{2}{r}j_n'(r) + \Bigl( 1 - \frac{n(n+1)}{r^2}\Bigr) j_n(r) = 0, \quad r > 0\,.
\end{equation*}
For every~$n \in \NN$ and every~$m \in \{-n, -n+1,\cdots, 0, \cdots, n\},$ the functions~$Y^m_n: \mathbb{S}^2 \to \RR$ are eigenfunctions of the Laplace-Beltrami operator on~$\mathbb{S}^2$ solving 
\begin{equation*}
	\Delta_{\mathbb{S}^2} Y^m_n + n(n+1) Y^m_n = 0 \quad \mbox{ on } \mathbb{S}^2\,.
\end{equation*}
Working in polar coordinates and writing~$x = r\omega, r = |x|, \omega \in \mathbb{S}^2,$ it follows by direct computation that~$u_{m,n} := j_n(kr)Y^m_n(\omega)$ solves Helmholtz equation~$\Delta u_{m,n} +k^2 u_{m,n} = 0.$ A general solution of Maxwell's equation in the unit ball~$B_1 \subset \R^3$ can be written down in terms of functions of the form~$u_{m,n}$ and its derivatives, and we will follow the treatment in~\cite[Chapter 2]{KH}. 

Define the vector fields 
\begin{equation*}
	U^m_n (\omega) = \frac{1}{\sqrt{n(n+1)}} \nabla_{\mathbb{S}^2} Y^m_n(\omega), \quad V^m_n(\omega) := \omega \times U^m_n(\omega)
\end{equation*}
These vector fields define an orthonormal frame of tangent fields on~$\mathbb{S}^2$ relative to the~$L^2(\mathbb{S}^2)$ inner product. Now suppose that~$f = \nu \times \bE$ on~$\mathbb{S}^2 = \partial B_1.$ Then it is shown in~\cite[Chapter 2]{KH} that the unique solution to Maxwell's equation with~$\e = 1$ and~$\mu = 1$ (and frequency~$k$) in~$D = B_1$ is given by 
\begin{multline} \label{e.maxwellE}
	\bE(r\omega) = \sum_{n=1}^\infty \sum_{m=-n}^n \bigl( f, V^m_n \bigr)_{L^2(\mathbb{S}^2)} \frac{\sqrt{n(n+1)}}{j_n(k) + k j_n'(k)} \frac{j_n(kr)}{r} Y^m_n(\omega) \omega \\
	\bigl(f, V^m_n\bigr)_{L^2(\mathbb{S}^2)} \frac{1}{j_n(k) + k j_n'(k)} \frac{j_n(kr) + k j_n'(kr)}{r} U^m_n(\omega)  - \bigl(f, U^m_n\bigr)_{L^2(\mathbb{S}^2)} \frac{j_n(kr)}{j_n(k)} V^m_n(\omega)\,. 
\end{multline}
and 
\begin{multline} \label{e.maxwellH}
	i k \bH(r\omega) = \sum_{n=1}^\infty \sum_{m=-n}^n \bigl( f, U^m_n\bigr)_{L^2(\mathbb{S}^2)}  \sqrt{n(n+1)} \frac{1}{j_n(k)} \frac{j_n(kr)}{r} Y^m_n(\omega )\omega \\
	+ \bigl(f, U^m_n \bigr)_{L^2(\mathbb{S}^2)}  \frac{1	}{j_n(k)} \frac{j_n(kr) + kj_n'(kr)}{r}U^m_n(\omega) - k^2 \bigl( f , V^m_n \bigr)_{L^2(\mathbb{S}^2)} \frac{j_n(kr)}{j_n(k) + k j_n'(k)} V^m_n(\omega)\,. 
\end{multline}
To be precise, the functions~$\bE,\bH$ defined above are the unique divergence free solution to the boundary value problem
\begin{equation*}
	\begin{dcases}
		\nabla \times \bE = - i k  \bH \quad \mbox{ in } D = B_1\\
		\nabla \times \bH = i k \bE \quad \mbox{ in } D = B_1\\
		\bE \times \nu_D = f\qquad \mbox{ on } \partial B_1\,.
	\end{dcases}
\end{equation*}
\begin{theorem} \label{t.electro}
	Let~$D = B_1 \subset \RR^3$, and let~$\Omega \supset \overline{D}$ denote an open bounded, Lipschitz domain. Let~$n \in \NN$ and~$m \in \Z \cap [-n,n].$ Then, setting~$k := k_n$ to be a zero of the Bessel function~$j_n,$ and define
	\begin{equation*}
		\begin{aligned}
			\bE_+(x) &:= \bE(r\omega) = \frac{\sqrt{n(n+1)}}{j_n(k_n) + k_n j_n'(k_n )} \frac{j_n(k_n r)}{r} Y^m_n(\omega) \omega \\
			&\qquad \qquad \qquad + \frac{1	}{j_n(k_n) + k_n j_n'(k_n)} \frac{j_n(k_n r) + k_n r j_n'(k_n r)}{r} U^m_n(\omega)\,,\\
			\bH_+(x) &:= ik_n \frac{j_n(k_nr)}{j_n(k_n) + k_n j_n'(k_n)} V^m_n(\omega)\,,
		\end{aligned}
	\end{equation*}
	for~$x \in D,$ so that the overdetermined set of boundary conditions~$\bE_+ \cdot \nu_D = 0$ at~$\partial D$ and~$\bH_+ = 0$ at~$\partial D$ hold. In addition, define~$(\bE_-,\bH_-)$ in~$\ENZ$ via 
	\begin{equation} \label{e.maxwellenz}
		\begin{dcases}
			\nabla \times \nabla \times \bE_- &= 0 \quad \mbox{ in } \ENZ\\
			\bE_- \times \nu_D &= \bE_+ \times \nu_D, \quad \mbox{ on } \partial D\\
			\bE_- \times \nu_\Omega &= 0 \quad \mbox{ on } \partial \Omega\,;
		\end{dcases}
	\end{equation} 
	There exists a unique solution to~\eqref{e.maxwellenz}  which is \emph{divergence free}, and which satisfies
	\begin{equation} \label{e.bcnormalization}
		\int_{\partial \Omega} \bE_- \cdot \nu_\Omega\,d\sh^2 = 0\,. 
	\end{equation}
	Finally, define
	\begin{equation*}
		\bH_- := \frac{1}{i k_n } \nabla \times \bE_- \quad \mbox{ in }\ENZ\,.
	\end{equation*}
	Then, setting 
	\begin{equation} \label{e.EHdef}
		\bE := \begin{dcases}
			\bE_+ \quad \mbox{ in } D\\
			\bE_- \quad \mbox{ in } \ENZ\,,
		\end{dcases}
		\quad \quad 
		\bH :=  \begin{dcases}
			\bH_+ \quad \mbox{ in } D\\
			0 \quad \mbox{ in } \ENZ\,,
		\end{dcases}
	\end{equation}
	defines an \emph{electrostatic} resonance of the Maxwell system (that is, solves~\eqref{e.vectorde=0}) with eigenvalue~$\la_0 := k_n^2.$  
\end{theorem}
\begin{proof}
	Fix~$n \in \NN,$ and let~$m \in \Z \cap [-n,n].$ Our construction of electrostatic resonances will be made by choosing~$f$ suitably so that the overdetermined boundary conditions at~$\partial B_1$ that lead to an electrostatic resonance are met. 
	To this end, set~$f := V^m_n,$ so that by the orthogonality of distinct spherical harmonics, the doubly infinite sums defining~$\bE$ and~$\bH$ collapse to relatively few terms each:
	\begin{equation} \label{e.inside}
		\bE(r\omega) = \frac{\sqrt{n(n+1)}}{j_n(k) + kj_n'(k)} \frac{j_n(kr)}{r} Y^m_n(\omega) \omega + \frac{1	}{j_n(k) + k j_n'(k)} \frac{j_n(kr) + kr j_n'(kr)}{r} U^m_n(\omega)\,,
	\end{equation}
	and 
	\begin{equation*}
		i k \bH(r\omega) = - k^2 \frac{j_n(kr)}{j_n(k) + kj_n'(k)} V^m_n(\omega)\,,
	\end{equation*}
	where, we will also make a special choice of~$k$ in what follows.
	
	At~$r=1$ we must arrange for~$\bE(\omega) \cdot\nu_D= 0$ and~$\bH = 0$. Both these goals are achieved upon selecting~$k$ to be a zero of~$j_n.$ The construction of the claimed electrostatic resonance in the inclusion~$D$ is complete. 
	
	Then, the existence and uniqueness of the solution to~\eqref{e.maxwellenz} that is divergence free and satisfies~\eqref{e.bcnormalization} now follows by a direct argument as follows: we seek~$\bE_- = \nabla h,$ so that the first equation in~\eqref{e.maxwellenz} holds. Then, to make~$\bE_-$ be divergence free, we ask that~$h$ be a harmonic function with suitable Dirichlet boundary conditions:
	\begin{equation*}
		\begin{dcases*}
			\Delta h= 0 \quad \mbox{ in } \ENZ\\
			\nabla h \times \nu_D= \bE_+ \times \nu_D = U^m_n (\omega) \times \nu_D  \quad \mbox{ on } \partial D\\
			h = \mbox{ constant }  \mbox{ on } \partial \Omega\,. 
		\end{dcases*}
	\end{equation*} 
	As~$\Omega$ and~$D$ are assumed to be simply connected, and since~$U^m_n = \frac{1}{\sqrt{n(n+1)}} \nabla_{\mathbb{S}^2} Y^m_n,$ there is only one free constant: we write 
	\begin{equation}
		h = h_0 + c_{\ENZ} \psi_{\ENZ}\,,
	\end{equation}
	with
	\begin{equation}
		\begin{dcases}
			\Delta h_0 = 0 \quad \mbox{ in } \ENZ\\
			h_0 = \frac{1}{\sqrt{n(n+1)}} Y^m_n \quad \mbox{ on } \partial D\\
			h_0= 0 \quad \mbox{ on } \partial \Omega\,,
		\end{dcases}
	\end{equation}
	the function~$\psi_{\ENZ}$ as identified via 
	\begin{equation}\label{e.psienz}
		\begin{dcases}
			\Delta \psi_{\ENZ} = 0 & \quad \mbox{ in } \ENZ\\
			\psi_{\ENZ} = 0 & \quad \mbox{ on } \partial D\\
			\psi_{\ENZ} = 1 & \quad \mbox{ on } \partial \Omega\,,
		\end{dcases}
	\end{equation} 
	and then the constant~$c_{\ENZ}$ chosen to satisfy~\eqref{e.bcnormalization}:
	\begin{equation*}
		c_{\ENZ} = - \frac{\int_{\partial \Omega} \partial_{\nu_\Omega} h_0 \,d\sh^2}{\int_{\partial \Omega} \partial_{\nu_\Omega} \psi_{\ENZ} \,d\sh^2}\,.
	\end{equation*}
	To see that this defines~$c_{\ENZ}$ to be a finite value, note that 
	since~$\psi_{\ENZ}$ is clearly a nonconstant function, we compute
	\begin{equation*}
		\int_{\partial \Omega} \frac{\partial \psi_{\ENZ}}{\partial \nu_\Omega} \,d \sh^2 = \int_{\partial (\ENZ)} \psi_{\ENZ} \frac{\partial \psi_{\ENZ}}{\partial \nu_\Omega} \,d \sh^2 = \int_{\ENZ} |\nabla \psi_{\ENZ}|^2 \neq 0\,. 
	\end{equation*} Since~$\bE_-$ is a gradient, it follows easily that~$\bH_- = 0$ in the~$\ENZ$ region. 
	
	It remains to check that~$\bE$ defined as in~\eqref{e.EHdef} satisfies~\eqref{e.strongform}. This is immediate, since~$\bE_+$ satisfies it by construction, with~$\la_0 := k_n^2,$ In addition since~$\bE_-$ is a gradient, it is curl-free, and so~\eqref{e.strongform} is also verified in~$\ENZ.$ Across~$\partial D,$  we constructed~$\bE_-$ so that it is tangent to~$\partial D,$ and that the tangential part of~$\bE$ is continuous across~$\partial D.$
	
	The proof of the Theorem is complete. 
\end{proof}
For completeness, here is the extension of the electrostatic resonance when~$\Omega$ is also spherical. If~$\Omega = B_R$ for some~$R > 1,$ then it can be checked by direct computation that defining~$\bE$ via~\eqref{e.inside} for~$r \leqslant 1$ and via~\eqref{e.outside} below for~$R \geqslant r > 1$ 
\begin{equation} \label{e.outside}
	\bE (r\omega) = \nabla h, \quad h (r,\omega) = \bigl( Ar^n + Br^{-n-1} \bigr) Y^m_n(\omega)\,,
\end{equation}
where~$A,B$ are chosen to satisfy~$A+B= \frac1{\sqrt{n(n+1)}}$ and~$AR^n + BR^{-n-1}=0.$ These last conditions enforce that~$\nabla_{\mathbb{S^2}} h = \frac1{\sqrt{n(n+1)}} \nabla_{\mathbb{S}^2} Y^m_n(\omega) = U^m_n(\omega)$ at~$r=1$, and~$h=0$ at~$r=R.$ Here,~$k$ is chosen to be a zero of~$j_n$, we obtain solution to~\eqref{e.strongform}. Note that we get~$\int_{\partial B_R} \bE \cdot \nu \,d\sh^2= \int_{\partial B_R} \nu \cdot \nabla_{\mathbb{S}^2} h \,d\sh^2 = 0, $ automatically, since the vector field~$U^m_n$ is tangent to~$\mathbb{S}^2.$ 

We close this section with a simple computation to show that even in the symmetrical case when~$D = B_1$ and~$\Omega = B_R,$ there exist eigenfunctions satisfying~\eqref{e.vectorde=0} that are \emph{not electrostatic}. 
\begin{theorem}
	\label{p.notelectrostatic}
	Let~$\Omega = B_R$ and~$D = B_1.$ Then there exist eigenfunctions satisfying~\eqref{e.vectorde=0} that are \emph{not electrostatic}, that is, the magnetic field~$\mathbf{H} = \nabla \times \bE$ does not identically vanish in~$\ENZ$. 
\end{theorem}
\begin{proof}
	As in the proof of the preceding theorem, we make a specific choice of~$f$ in~\eqref{e.maxwellE} and~\eqref{e.maxwellH}; however, contrary to that case, fixing~$p \in \mathbb{N}$ and~$q \in \Z \cap [-p,p],$ we set~$f = U^q_p$ (instead of the choice~$f = V^q_p$ made in the proof of the previous theorem). With this choice, the solutions to Maxwell's equation in the inclusion~$D$ yields 
	\begin{equation} \label{e.insideD}
		\begin{dcases}
			\bE_+ (r\omega) &= - \frac{j_p(kr)}{j_p(k)} V^q_p(\omega), \quad r \leqslant 1\\
			ik \bH_+(r\omega) &= \sqrt{p(p+1)} \frac{1}{j_p(k)} \frac{j_p(kr)}{r} Y^q_m(\omega)\omega + \frac{1}{j_p(k)} \frac{j_p(kr) + k j_p'(kr)}{r} U^q_p(\omega)\,, \quad r \leqslant 1\,. 
		\end{dcases}
	\end{equation}
	With this choice, we have arranged that~$\bE_+ \cdot \nu_D = 0,$ and we must extend~$(\bE,\bH)$ to~$\ENZ$ so that the tangential of both~$\bE$ and~$\bH$ are continuous across~$\partial D.$ This requirement will select~$k,$ which has so far not been pinned down. 
	
	In order to extend~$\bE = \bE_-$ in the ENZ region~$\ENZ$, a natural idea is to solve the boundary value problem~\eqref{e.maxwellenz} with data coming from~$D:$
	\begin{equation*}
		\begin{dcases}
			\nabla  \times \bE_-= i k \bH_-\quad \mbox{ in } \ENZ = B_R \setminus \overline{B_1}\\
			\nabla \times \bH_- = 0 \qquad \qquad \mbox{ in } \ENZ= B_R \setminus \overline{B_1}\\
			\nu_D \times \bE_- = \omega \times \bE_- = - \omega \times V^q_p(\omega) = -\frac{1}{\sqrt{p(p+1)}} \omega \times (\omega \times U^q_p(\omega)) \\
			\qquad \qquad \qquad \qquad \qquad \qquad \qquad \ \quad = \frac{1}{\sqrt{p(p+1)}} U^q_p(\omega)\, \quad \mbox{ on } r = 1\\
			\nu_\Omega  \times \bE_- = 0 \quad \mbox{ on } r = R\\
			\mathrm{div}\, \bE_- = 0 \quad \mbox{ in } B_R\setminus \overline{B_1}\,.
		\end{dcases}
	\end{equation*}
	
	We solve this problem by separation of variables. By~\cite[Corollary 2.47]{KH} we are led to look for a solution (which will turn out to be the unique one) of the form
	\begin{equation*}
		\bE_-(r,\omega) = a^q_p(r) Y^q_p(\omega)\omega + b^q_p(r) U^q_p(\omega) + c^q_p(r) V^q_p(\omega)\,,
	\end{equation*}
	with 
	\begin{equation*}
		a^q_p(r) := \int_{\mathbb{S}^2} \omega \cdot \bE(r\omega) Y^{-q}_p(\omega) \,d \sh^2(\omega) = \frac1{r} \int_{\mathbb{S}^2} x\cdot \bE(x) Y^{-q}_p(\omega)\,d\sh^2(\omega)\,.
	\end{equation*}
	Therefore,~$ra^q_p(r)$ is the coefficient of the function~$x\mapsto x\cdot \bE(x)$ at mode~$Y^q_p(\omega).$
	
	We now make an observation about this function. Since~$\bE$ is divergence-free and~$0 = \nabla \times \nabla \times \bE = - \Delta \bE=0,$ we claim that~$x \cdot \bE(x)$ is harmonic (as a real-valued function). To see this, we compute
	\begin{equation*}
		\Delta (x \cdot \bE) = \sum_{j=1}^3\sum_{i=1}^3 \partial_j^2 x_i \bE_i = \sum_{j=1}^3\sum_{i=1}^3 \partial_j \bigl( \delta_{ij} \bE_i + x_i \partial_j \bE_i\bigr) = \sum_{j=1}^3 \sum_{i=1}^3 \delta_{ij} \partial_j\bE_i + x_i \partial_j^2 \bE_i = 0\,.
	\end{equation*}
	It follows by separation of variables (see~\cite[Theorem 2.22, for e.g.]{KH}) that
	\begin{equation} \label{e.aqp}
		ra^q_p(r) = \bigl( A^q_p r^p + B^q_p r^{-p-1}\bigr) \,,
	\end{equation}
	for some~$A^q_p,B^q_p \in \RR.$
	
	Next, we derive conditions on the coefficients~$a^q_p(r), b^q_p(r)$ and~$c^q_p(r)$ that ensure that~$\mathrm{div}\, \bE_-=0.$ In polar coordinates, we obtain
	\begin{equation*}
		\mathrm{div}\, \bE = \frac1{r} \partial_r \bigl( r^2 \bE_r\bigr) + \mathrm{div}_{\mathbb{S}^2} \, \bE_{tan} = 0\,, 
	\end{equation*}
	and so, 
	\begin{align*}
		\frac{1}{r}\partial_r \bigl( r^2 \bE_r\bigr) &= \int_{\mathbb{S}^2} \partial_r \bigl( r^2 \bE_r\bigr) Y^{-q}_p(\omega) \,d\sh^2(\omega) \\
		& = - \int_{\mathbb{S}^2} \mathrm{div}_{\mathrm{S}^2}\, \bE_{tan} (r,\omega) \,d\sh^2(\omega) \\
		&= \int_{\mathrm{S}^2} \bE(r,\omega) \cdot \nabla_{\mathbb{S}^2} Y^{-q}_p(\omega)\,d\sh^2(\omega) \\
		& = \int_{\mathrm{S}^2} \bE(r,\omega) \cdot \sqrt{p(p+1)} U^{-q}_p(\omega)\,d\sh^2(\omega) \\
		& =  \sqrt{p(p+1)} b^q_p(r)
	\end{align*}
	Inserting this in~\eqref{e.aqp}, we arrive at 
	\begin{equation} \label{e.bqp}
		b^q_p(r) = \frac{1}{\sqrt{p(p+1)}} \partial_r\bigl(r^2 a^q_p(r)\bigr) = \sqrt{\frac{p+1}{p}} A^q_p r^{p-1} - \sqrt{\frac{p}{p+1}} B^q_p r^{-p-2}\,. 
	\end{equation}
	Finally, for the coefficients~$c^q_p(r)$ we obtain
	\begin{align*}
		\sqrt{p(p+1)}c^q_p(r) &= \int_{\mathbb{S}^2} \bE(r,\omega) \cdot( \omega \times \nabla_{\mathbb{S}^2} Y^{-q}_p(\omega)) \,d\sh^2(\omega)\\
		&= \int_{\mathbb{S}^2} \bE(r,\omega) \times \omega \cdot \nabla_{\mathbb{S}^2} Y^{-q}_p(\omega)\,d\sh^2(\omega)\\
		& = - \int_{\mathbb{S}^2} \mathrm{div}_{\mathbb{S}^2} \bigl( \bE(r\omega) \times \omega\bigr) Y^{-q}_p(\omega)\,d\sh^2(\omega) \\
		&= - r \int_{\mathbb{S}^2} \omega \cdot \mathrm{curl} \,  \bE(r\omega) Y^{-q}_p(\omega) \,d\sh^2(\omega) \\
		&= - \int_{\mathbb{S}^2} x \cdot \mathrm{curl} \, \bE(r\omega) Y^{-q}_p(\omega)\,d\sh^2(\omega)
	\end{align*}
	where, in the last two lines we used the vector calculus identity from~\cite[Corollary 6.20]{KH}:
	\begin{equation*}
		\mathrm{div}_{\mathbb{S}^2} ( \bE (r,\omega) \times \omega) = \mathrm{div}_{\mathbb{S}^2}|_{|x| = r} \bigl( \bE(r\omega) \times \omega \bigr) = r \omega \cdot \mathrm{curl}\, \bE(r\omega)\,. 
	\end{equation*}
	This computation implies that~$\sqrt{p(p+1)} c^q_p(r)$ is the expansion coefficient of the scalar function~$\psi(x) := x \cdot \mathrm{curl}\, \bE(x).$ Since~$\bE$ is~$\mathrm{curl}\, \mathrm{curl}-$free and divergence free, it follows that~$\psi$ is harmonic, by a similar computation as above. This implies that there exist constants~$C^q_p,D^q_p \in \RR$ such that 
	\begin{equation*}
		\sqrt{p(p+1)} c^q_p(r) = C^q_p r^p + D^q_p r^{-p-1}\,. 
	\end{equation*}
	
	Putting this all together, we conclude that 
	\begin{multline} \label{e.E-}
		\bE_- (r,\omega) = \bigl(A^q_p r^{p-1} + B^q_p r^{-p-2} \bigr) Y^q_p(\omega)\omega + \Bigl( \sqrt{\frac{p+1}{p}} A^q_p r^{p-1} - \sqrt{\frac{p}{p+1}} B^q_p r^{-p-2}\Bigr) U^q_p(\omega) \\
		+ \frac{1}{\sqrt{p(p+1)}} \bigl( C^q_p r^p + D^q_p r^{-p-1} \bigr) V^q_p(\omega)\,. 
	\end{multline}
	It remains to impose boundary conditions and solve for~$A^q_p, B^q_p, C^q_p, D^q_p,$ and then impose continuity of the magnetic field~$\bH$ to select~$k.$ 
	
	\smallskip 
	Since~$\bE_+(r,\omega) = -\frac{j_p(kr)}{j_p(k)} V^q_p(\omega)$ for~$r \leqslant 1,$ continuity of the tangential component of~$\bE$ at~$r=1$ requires
	\begin{equation*}
		\begin{dcases}
			-1 &=  \frac{1}{\sqrt{p(p+1)}} (C^q_p + D^q_p) \qquad \qquad  \mbox{[Continuity of the~$V^q_p$ component]}\\
			0 &= \sqrt{\frac{p+1}{p}} A^q_p - \sqrt{\frac{p}{p+1}} B^q_p \qquad \qquad  \mbox{[Continuity of the~$U^q_p$ component]}\,.
		\end{dcases} 
	\end{equation*}
	Next, imposing the perfect electrical conducting boundary conditions at~$r=R$ we obtain the following system of equations:
	\begin{equation*}
		\begin{dcases}
			0 &= \frac{\sqrt{p+1}}{\sqrt{p}} A^q_p R^{p-1} - \sqrt{\frac{p}{p+1}} R^{-p-2} \\
			0 &=  C^q_p R^p + D^q_p R^{-p-1}\,. 
		\end{dcases}
	\end{equation*}
	Evidently~$A^q_p = B^q_p =0,$ and~$C^q_p,D^q_p$ are determined by 
	\begin{equation*}
		C^q_p + D^q_p = - \sqrt{p(p+1)}, \quad C^q_p R^{p} + D^q_p R^{-p-1} = 0\,. 
	\end{equation*}
	This linear system can be uniquely inverted to obtain~$C^q_p,D^q_p,$ as a function of~$R$, and this solution does not depend on~$k,$ and we obtain 
	\begin{equation*}
		\bE_-(r,\omega) = \frac{1}{\sqrt{p(p+1)}} \bigl( C^q_p r^p + D^q_p r^{-p-1} \bigr) V^q_p(\omega)\,.
	\end{equation*}
	It remains to check that~$\bE_-$ satisfies~$\nabla \times \nabla \times \bE_-(r,\omega) = 0$ (it is automatically divergence free by construction). It is possible to check this following~\cite[Theorem 2.48]{KH}, or more directly, through a tedious computation.

Having solved for these unknown coefficients, we can use Maxwell's equations to determine that the magnetic field is given by
\begin{multline*}
	\bH_- (r,\omega) = - \frac{1}{i k} \Bigl\{   \bigl( C^q_p r^{p-1} + D^q_p r^{-p-2} \bigr) Y^q_p (\omega) \omega \\+   \Bigl( (p+1)C^q_p r^{p-1} - pD^q_p r^{-p-2}) \Bigr) U^q_p(\omega)\Bigr\} \quad 1 \leqslant r \leqslant R\,.  
\end{multline*}
Recalling from~\eqref{e.insideD} that~$\bH_+$ is given by 
\begin{multline*}
	\bH_+(r,\omega) = \frac1{ik} \biggl[\sqrt{p(p+1)} \frac{1}{j_p(k)} \frac{j_p(kr)}{r} Y^q_m(\omega)\omega \\+ \frac{1}{j_p(k)} \frac{j_p(kr) + k j_p'(kr)}{r} U^q_p(\omega) \biggr]\qquad \qquad 0 < r \leqslant 1\,,
\end{multline*}
continuity of the tangential component of~$\bH$ across~$r=1$ needs 
\begin{equation} \label{e.ivt}
	-\bigl((p+1) C^q_p - D^q_p\bigr) = \frac{j_p(k) + k j_p'(k)}{j_p(k)} = 1 + k \frac{j_p'(k)}{j_p(k)}\,.
\end{equation}
On the other hand, continuity of the normal component of~$\bH$ is automatic since~$C^q_p + D^q_p = - \sqrt{p(p+1)}.$ 

As the roots of~$j_k$ and~$j_k'$ interlace strictly, and as the continuous function~$k \mapsto k \frac{j_p'(k)}{j_p(k)}$ mapping onto~$(-\infty,\infty)$ as~$k$ varies between any pair of zeroes of the Bessel function~$j_p(k),$ the intermediate value theorem implies that there is a solution to the equation~\eqref{e.ivt}, for any fixed~$R,p$. This is because, fixing~$R,p,$ uniquely determines the quantity~$-\bigl( (p+1) C^q_p - D^q_p\bigr) - 1$ independently of~$k.$ 

This concludes the proof of the Theorem. 
\end{proof}

	\appendix
\section{Rellich's theorem} \label{sec:rellich}
We repeatedly appeal to a powerful theorem of Rellich~\cite[Theorem 1 of Chapter II, Sec. 2]{Rel}, whose statement we recall for the reader's convenience. 

\begin{theorem} \label{t.rellich}
	Suppose that~$\{A(\de)\}_{|\de| < \de_0 }$ bounded operator on a Hilbert space~$H$ which is a convergent power series in~$\de$ defined in some neighborhood~$\{\de \in \CC: |\de| < \de_0\} \subset \CC,$ given by 
	\begin{equation*}
		A(\de) = A_0 + \de A_1 + \cdots \,. 
	\end{equation*} 
	Suppose that for real~$\de$ with~$|\de|$ small and for any~$u,v \in H$ we have 
	\begin{equation} \label{e.selfadj}
		\bigl( A(\de) u, v \bigr) = \bigl(u, A(\de) v\bigr)\,,
	\end{equation}
	with~$(\cdot, \cdot)$ denoting the inner product in~$H$. Suppose that~$\la$ is an eigenvalue of finite multiplicity~$h$ of the operator~$A(0) =: A_0$, and suppose there are positive numbers~$d_1, d_2$ such that the spectrum of~$A_0$ in the open interval~$(\la_0 - d_1, \la_0 + d_2)$ consists exactly of the point~$\la.$ Then there exists~$h$ families of complex numbers 
	\begin{equation*}
		\la_1(\de), \cdots, \la_h(\de)\,,
	\end{equation*}
	that have convergent power series in~$\de$, and associated elements~$\{\phi_1(\de),\cdots, \phi_h(\de)\} \subset H$ also given by a convergent power series in the Hilbert space~$H$ that satisfy the following properties: 
	\begin{enumerate}
		\item The element~$\phi_i(\de)$ is an eigenvector of~$A(\de)$ corresponding to~$\la_i(\de)$ as eigenvalue:
		\begin{equation*}
			A(\de) \phi_i(\de) = \la_i(\de)\phi_i(\de)\,, \quad i = 1, \cdots, h\,.
		\end{equation*}
		Furthermore, for each~$i \in \{1,\cdots, h\}$, we have that~$\la_i(0) = \la$ and in addition, for all~$\de \in \R$ we have that 
		\begin{equation*}
			(\phi_i(\de), \phi_j(\de)) = \de_{ij} = \begin{dcases}
				1 & \mbox{ if } i = j\\
				0 & \mbox{ otherwise } 
			\end{dcases}\,. 
		\end{equation*}
		\item For each pair of positive numbers~$d_1', d_2'$ with~$d_1' < d_1, d_2' < d_2$ there exists~$\rho > 0$ such that the spectrum of~$A(\de)$ in~$\la - d_1'\leqslant \mu \leqslant \la + d_2'$ for all real~$\de$ with~$|\de| < \rho$ consists exactly of the points~$\la_1(\de), \cdots, \la_h(\de).$ 
	\end{enumerate}
\end{theorem}

	\bibliographystyle{alpha}
\bibliography{ref.bib}
\end{document}